\documentclass[A4paper,12pt]{book}
\usepackage{amssymb,amsmath}
\usepackage[russian]{babel}
\usepackage{graphicx}
\newcount\exnom\exnom=0

\long\def\exo/{\vspace{0.2cm} \noindent\advance\exnom by1{\bf
{\the\exnom}.}}

\newcommand{\ds}{\displaystyle}

\begin{document}

\begin{center}\large  \textbf{DERIVATIONS IN SOME FINITE ENDOMORPHISM SEMIRINGS}
\end{center}

\begin{center} IVAN TRENDAFILOV 
\end{center}

\thispagestyle{empty}

\begin{quote} Abstract.{\small
 The goal of this paper is to provide some basic structure information on derivations in finite semirings.}
\end{quote}

\vspace{5mm}

\centerline{{\large 1. \hspace{0.5mm}Introduction}}

\vspace{3mm}

Over a period of sixteen years differential algebra
went from being an approach that many people mistrusted or misunderstood
to being a part of algebra that enjoys almost
unquestioned acceptance. This algebra has been studied by many authors for the last 60 years and especially the relationships between
derivations and the structure of rings. The notion of the ring with derivation is quite old and plays an important role in the integration of
analysis, algebraic geometry and algebra. In the 1940's it was discovered that the Galois theory of algebraic
equations can be transferred to the theory of ordinary linear differential equations (the Picard-Vessiot theory).
In the 1950's the differential algebra was initiated by the works of J. F. Ritt and E. R. Kolchin.
In 1950, Ritt [11] and in 1973, Kolchin [9] wrote the  classical  books on differential algebra.

The theory of derivations plays a significant role not
only in ring theory, but also in functional analysis and  linear differential equations.  For instance, the classical Noether-Skolem theorem yields the solution of the problem for finite dimensional central simple algebras (see the well-known Herstein's book [6]).
One of the natural questions  in algebra and analysis is whether a map can be
defined by its local properties. For example, the question whether a map which acts like a derivation on
the Lie product of some important Lie subalgebra of prime rings is induced by an ordinary derivation was
a well-known problem posed by Herstein [5]. The first result in this direction was obtained in an unpublished
work by Kaplansky, who considered matrix algebras over a field.
 Herstein's problem was solved in full generality only after the powerful technique of functional identities was developed, see [1].
In 1950's, Herstein ([4], [5], [7]) started the study of the relationship between the associative structure and the Jordan and Lie
structures of associative rings.

An additive mapping $D$ from $R$ to $R$, where $R$ is an associative ring, is called a Jordan derivation if $D(x^2) = D(x)x + xD(x)$ holds for all $x \in R$. Every derivation is obviously a Jordan derivation and the converse is in general not true.

It is important to note that the definition of Jordan derivation presented in the work by Herstein is not the same as the one given above. In fact, Herstein constructed, starting from the ring $R$, a new ring, namely the Jordan ring R, defining the product in this $a\circ b = ab + ba$ for any $a, b \in R$. This new product is well-defined and it can be easily verified that $(R,+,\circ)$ is a ring. So, an additive mapping $D$, from
the Jordan ring into itself, is said by Herstein to be a Jordan derivation, if $D(a\circ b) = D(a)\circ b + a\circ D(b)$, for
any $a, b \in R$. So, in 1957, Herstein proved a classical result: "If $R$ is a prime ring of a characteristic different from 2, then every
Jordan derivation of $R$ is a derivation."

During the last few decades there has been a great deal of works concerning derivations $D_i$ in rings, in Lie rings, in skew polynomial rings and other structures, which commutes, i.e. $D_iD_j = D_jD_i$.

What we know about derivatins in semirings? Nothing, or almost nothing exept the definition in Golan's book [3] and a few propositions.

This paper is an attempt to start a study of derivations in finite semirings. Following Herstein's idea of multiplication in Jordan ring, we construct derivations in the endomorphism semiring of a finite chain.

The paper is organized as follows.
After the second section of preliminaries, in section 3 we introduce a semiring consisting of endomorphisms having an image with two fixed elements called a string. In such a string we consider the arithmetic and some kinds of nilpotent elements and subsemirings. In section 4 we construct a mapping $D$ from the given string into itself and prove that $D$ is a derivation in one subsemiring of the string. Then we show that the semiring is a maximal differential subsemiring of this string.
Section 5 is devoted to the construction of  maps $\delta_\alpha$ from given string into itself. They are Jordan multiplications, and we are studying their properties. The main results are that $\delta_\alpha$ are derivations which commute and that the set of all derivations is a multiplicative semilattice with an identity and an absorbing element. In section 6 we generalize the notion of string and consider the arithmetic in such strings. In section 7 we give some counterexamples and show that the maps $\delta_\alpha$, where $\alpha$ are from the whole string, are derivations in an ideal of this string. Finally in this section we consider a class of maps which are the derivations in the whole string.

\vspace{5mm}

\centerline{\large 2. \hspace{0.5mm}Preliminaries}

\vspace{3mm}

 An algebra $R = (R,+,.)$ with two binary operations $+$ and $\cdot$ on $R$, is called {semiring} if:

$\bullet\; (R,+)$ is a commutative semigroup,

$\bullet\; (R,\cdot)$ is a semigroup,

$\bullet\;$ both distributive laws hold
$ x\cdot(y + z) = x\cdot y + x\cdot z$ and $(x + y)\cdot z = x\cdot z + y\cdot z$
for any $x, y, z \in R$.

 Let $R = (R,+,.)$ be a semiring.
If a neutral element $0$ of  semigroup $(R,+)$ exists and satisfies $0\cdot x = x\cdot 0 = 0$ for all $x \in R$, then it is called zero.
If a neutral element $1$ of  semigroup $(R,\cdot)$ exists, it is called one.

An element $a$ of a semiring $R$ is called additively (multiplicatively) idempotent if ${a + a = a}$ $\;(a\cdot a = a)$.
A semiring $R$ is called additively idempotent if each of its elements is additively idempotent.

An element $a$ of a semiring $R$ is called an additively (multiplicatively) absorbing element if and only if
$a + x = a\;\;$ $(a\cdot x = x\cdot a = a)\;\;$ for any $x \in R$.
The zero of $R$ is the unique multiplicative absorbing element; of course it does not need to exist.
Following [10], an element  of a semiring $R$
is called an infinity if it is both additively and multiplicatively absorbing.

\vspace{2mm}

Facts concerning semirings, congruence relations in semirings and (right, left) ideals of semirings  can be found in [3] and [10].

\vspace{2mm}
 An algebra$\mathcal{M}$ with binary operation $\vee$ such as

$\bullet$ $\; a\vee(b\vee c) = (a\vee b)\vee c$ for any $a, b, c \in \mathcal{M}$;

$\bullet$ $\; a\vee b = b \vee a$ for any $a, b\in \mathcal{M}$;

$\bullet$ $\; a\vee a = a$ for any $a \in \mathcal{M}$.

is called semilattice (join semilattice).

Another term used for $\mathcal{M}$ is a commutative idempotent semigroup -- see [15].
For any $a, b\in \mathcal{M}$ we denote
$ a \leq b \; \iff \; a \vee b = b$.
In this notation, if there is a neutral element in $\mathcal{M}$, it is the least element.

\vspace{2mm}

For a semilattice $\mathcal{M}$ the set $\mathcal{E}_\mathcal{M}$ of the endomorphisms of $\mathcal{M}$ is a semiring
 with respect to the addition and multiplication defined by:

 $\bullet \; h = f + g \; \mbox{when} \; h(x) = f(x)\vee g(x) \; \mbox{for all} \; x \in \mathcal{M}$,

 $\bullet \; h = f\cdot g \; \mbox{when} \; h(x) = f\left(g(x)\right) \; \mbox{for all} \; x \in \mathcal{M}$.

 This semiring is called the \textbf{\emph{ endomorphism semiring}} of $\mathcal{M}$.
It is important to note that in this paper all semilattices are finite chains.
Following [12] and [13] we fix a finite chain $\mathcal{C}_n = \; \left(\{0, 1, \ldots, n - 1\}\,,\,\vee\right)\;$ and denote the endomorphism semiring of this chain with $\widehat{\mathcal{E}}_{\mathcal{C}_n}$. We do not assume that $\alpha(0) = 0$ for arbitrary $\alpha \in \widehat{\mathcal{E}}_{\mathcal{C}_n}$. So, there is not a zero in  endomorphism semiring $\widehat{\mathcal{E}}_{\mathcal{C}_n}$. Subsemiring ${\mathcal{E}}_{\mathcal{C}_n} = {\mathcal{E}}_{\mathcal{C}_n}^0$ of $\widehat{\mathcal{E}}_{\mathcal{C}_n}$ consisting of all maps $\alpha$ with property $\alpha(0) = 0$ has zero and  is considered in [12] and [15].

\vspace{2mm}

If $\alpha \in \widehat{\mathcal{E}}_{\mathcal{C}_n}$ such that $f(k) = i_k$ for any  $k \in \mathcal{C}_n$ we denote $\alpha$ as an ordered $n$--tuple $\wr\,i_0,i_1,i_2, \ldots, i_{n-1}\,\wr$. Note that mappings will be composed
accordingly, although we shall usually give preference to writing mappings on
the right, so that $\alpha \cdot \beta$ means "first $\alpha$, then $\beta$".

For other properties of the endomorphism semiring we refer to [8], [12],  [13] and [15].

In the following sections we use some terms from book [2] having in mind that in [14] we show that some subsemigroups of the partial transformation semigroup are indeed endomorphism semirings.

\vspace{5mm}

\centerline{\large 3. \hspace{0.5mm} Strings of type 2}

\vspace{3mm}

By $\mathcal{STR}\{a,b\}$ we denote subset of $\widehat{\mathcal{E}}_{\mathcal{C}_n}$ consisting of endomorphisms with image $\{a,b\}$ (either $\{a\}$, or $\{b\}$) where $a, b \in \mathcal{C}_n$. This set is called a \textbf{\emph{string of type 2}}. So, for fixed $a, b \in \mathcal{C}_n$ and $a < b$, a string of type 2 is
$$\mathcal{STR}\{a,b\} = \{\,\wr\,a,\ldots,a\,\wr, \wr\,a,\ldots,a,b\,\wr, \ldots, \wr\,a,b,\ldots,b\,\wr, \wr\,b,\ldots,b\,\wr\,\}.$$

In the semiring $\widehat{\mathcal{E}}_{\mathcal{C}_n}$ there is an order of the following way:

For two arbitrary endomorphisms $\alpha = \wr\,k_0,k_1,\ldots,k_{n-1}\,\wr$ and $\beta = \wr\,l_0,l_1,\ldots,l_{n-1}\,\wr$ the relation $\alpha \leq \beta$ means that $k_i \leq l_i$ for all $i = 0,1,\ldots,n-1$, i.e. $\alpha + \beta = \beta$.

 With regard to this order each of  sets $\mathcal{STR}\{a,b\}$ is an $n+1$ -- element chain with the least element $\wr\,a,\ldots,a\,\wr$ and the biggest element $\wr\,b,\ldots,b\,\wr$. Hence string $\mathcal{STR}\{a,b\}$ is closed under the addition of  semiring $\widehat{\mathcal{E}}_{\mathcal{C}_n}$. It is easy to see that the composition of two endomorphisms with images $\{a,b\}$ is an endomorphism of such type. So, we have

\vspace{3mm}

\textbf{Proposition 3.1 } \textsl{For any $a, b \in \mathcal{C}_n$  string $\mathcal{STR}\{a,b\}$ is a subsemiring of  semiring $\widehat{\mathcal{E}}_{{\mathcal{C}_n}}$.}

\vspace{3mm}

Note that  string $\mathcal{STR}\{a,b\}$ is not a subsemiring of semiring ${\mathcal{E}}_{\mathcal{C}_n}^{(a)}\cap {\mathcal{E}}_{\mathcal{C}_n}^{(b)}$ (see [13]), namely $\wr\,2,2,3,3\,\wr \in {\mathcal{E}}_{\mathcal{C}_4}^{(3)}$, but $\wr\,2,2,3,3\,\wr \notin {\mathcal{E}}_{\mathcal{C}_4}^{(2)}$.

\vspace{3mm}

Let us denote the elements of  semiring $\mathcal{STR}\{a,b\}$ by $\alpha^{a,b}_k$, where $k = 0,\ldots,n$ is the number of the elements of $\mathcal{C}_n$ with an image equal to $b$, i.e.
$$\alpha^{a,b}_k = \wr\, a, \ldots, a, \underbrace{b, \ldots, b}\,\wr.$$

\vspace{-3mm}
{\scriptsize
\centerline{$\hphantom{aaaaaaaaaaaaaaaaaaa} k \hphantom{aa} $}
}

When $a$ and $b$ are fixed, we replace $\alpha^{a,b}_k$ by $\alpha_k$.

\vspace{3mm}

\textbf{Proposition 3.2 } \textsl{Let $a, b \in \mathcal{C}_n$ and $\mathcal{STR}\{a,b\} = \{\alpha_0, \ldots, \alpha_n\}$. For every $k = 0,\ldots,n$ follows}
$$\begin{array}{lll}
\alpha_k \cdot \alpha_s = \alpha_0, & \mbox{\textsl{where}} & 0 \leq s \leq n - b - 1\\
\alpha_k \cdot \alpha_s = \alpha_k, & \mbox{\textsl{where}} & n-b \leq s \leq n - a - 1\\
\alpha_k \cdot \alpha_s = \alpha_n, & \mbox{\textsl{where}} & n-a \leq s \leq n
\end{array}.$$

\emph{Proof.} Let us multiply $\alpha_k \cdot \alpha_s$, where $0 \leq s \leq n - b - 1$. This means that at least $b+1$ elements of $\mathcal{C}_n$ have images under endomorphism $\alpha_s$ equal to $a$ and then $b$ is not a fixed point of $\alpha_s$. Hence, for all $i \in \mathcal{C}_n$ follows
$$(\alpha_k \cdot \alpha_s)(i) = \alpha_s(\alpha_k(i)) = \left\{ \begin{array}{ll} \alpha_s(a),& \mbox{if}\; i \leq n -k - 1\\ \alpha_s(b), & \mbox{if} \; i \geq n - k \end{array} = a \right.$$
which means that $\alpha_k \cdot \alpha_s = \alpha_0$.

Let us multiply $\alpha_k \cdot \alpha_s$, where $n - b \leq s \leq n - a - 1$. Now $a$ and $b$ are both fixed points of the endomorphism $\alpha_s$. Hence, for all $i \in \mathcal{C}_n$ follows
$$(\alpha_k \cdot \alpha_s)(i) = \alpha_s(\alpha_k(i)) = \left\{ \begin{array}{ll} \alpha_s(a),& \mbox{if}\; i \leq n -k-1\\ \alpha_s(b), & \mbox{if} \; i \geq n - k \end{array} = \left\{ \begin{array}{ll} a,& \mbox{if}\; i \leq n -k-1\\ b, & \mbox{if} \; i \geq n - k \end{array} = \alpha_k(i) \right. \right.$$
which means that $\alpha_k \cdot \alpha_s = \alpha_k$.

 Let us multiply $\alpha_k \cdot \alpha_s$, where $n - a \leq s \leq n$. This means that at least $n - a$ elements of $\mathcal{C}_n$ have images under endomorphism $\alpha_s$ equal to $b$ and then $a$ is not a fixed point of $\alpha_s$. Hence, for all $i \in \mathcal{C}_n$ follows
$$(\alpha_k \cdot \alpha_s)(i) = \alpha_s(\alpha_k(i)) = \left\{ \begin{array}{ll} \alpha_s(a),& \mbox{if}\; i \leq n -k-1\\ \alpha_s(b), & \mbox{if} \; i \geq n - k \end{array} = b \right.$$
which means that $\alpha_k \cdot \alpha_s = \alpha_n$.

\vspace{3mm}

The endomorphism $\alpha$ is called $\mathbf{a}$ \textbf{-- \emph{nilpotent}} if for some natural $k$ follows that $\alpha^k = \alpha_0$ and respectively, $\mathbf{b}$ \textbf{-- \emph{ nilpotent}} if $\alpha^k = \alpha_n$. From the last proposition follows that in both cases $k = 2$ and

\vspace{3mm}

\textbf{Corollary 3.3 } \textsl{For any $a, b \in \mathcal{C}_n$ the set:}

\textsl{a.  $N_a = \{\alpha_0, \ldots, \alpha_{n-b-1}\}$ of all $a$ - nilpotent elements is a subsemiring  of  semiring $\mathcal{STR}\{a,b\}$.}

\textsl{b.  $N_b = \{\alpha_{n-a}, \ldots, \alpha_{n}\}$ of all $b$ - nilpotent elements is a subsemiring  of  semiring $\mathcal{STR}\{a,b\}$.}

\textsl{c.  $I\!d_{a,b} = \{\alpha_{n-b}, \ldots, \alpha_{n-a-1}\}$ of all idempotent elements is a subsemiring  of  semiring $\mathcal{STR}\{a,b\}$.}

\vspace{3mm}

Note that all idempotent elements in  semiring $\widehat{\mathcal{E}}_{\mathcal{C}_n}$ do not form a semiring.

\vspace{3mm}

The constant endomorphisms $\wr\, a, \ldots, a \,\wr$ and $\wr\, b, \ldots, b \,\wr$ are called \textbf{\emph{centers}} of  semirings ${\mathcal{E}}_{\mathcal{C}_n}^{(a)}$ and ${\mathcal{E}}_{\mathcal{C}_n}^{(b)}$, respectively.
 So, we may imagine that  string
$\mathcal{STR}\{a,b\}$ "connects" the two centers $\wr \, a, \ldots, a \,\wr$ and $\wr\, b, \ldots, b \,\wr$.
 Semiring $N_a$ is a subsemiring of ${\mathcal{E}}_{\mathcal{C}_n}^{(a)}$, analogously $N_b$ is a subsemiring of ${\mathcal{E}}_{\mathcal{C}_n}^{(b)}$ and semiring $I\!d_{a,b}$ is a subsemiring of ${\mathcal{E}}_{\mathcal{C}_n}^{(a)}\cap {\mathcal{E}}_{\mathcal{C}_n}^{(b)}$.

\vspace{3mm}

Let us consider the subset $S_{a,b} = \{\alpha_0, \ldots, \alpha_{n-a-1}\}$ of the  $\mathcal{STR}\{a,b\}$. Obviously, the set $S_{a,b}$ is closed under the addition. So, immediately from Proposition 3.2 follows

\vspace{3mm}

\textbf{Corollary 3.4 } \textsl{For any $a, b \in \mathcal{C}_n$ the set $S_{a,b}$ is a subsemiring of  semiring $\mathcal{STR}\{a,b\}$.}

\vspace{3mm}

Note that  semiring $S_{a,b}$ is a disjoint union of  semirings $N_a$ and $I\!d_{a,b}$.

\vspace{3mm}

From dual point of view we consider the subset $T_{a,b} = \{\alpha_{n-b}, \ldots, \alpha_{n}\}$ of  string $\mathcal{STR}\{a,b\}$. This set is also closed under the addition and from Proposition 3.2 we have

\vspace{3mm}

\textbf{Corollary 3.5 } \textsl{For any $a, b \in \mathcal{C}_n$ the set $T_{a,b}$ is a subsemiring of  semiring $\mathcal{STR}\{a,b\}$.}

\vspace{3mm}

Note that  semiring $T_{a,b}$ is a disjoint union of  semirings  $I\!d_{a,b}$ and $N_b$.

\vspace{5mm}

 \centerline{\large 4. \hspace{0.5mm}Derivation ${D}$}

\vspace{3mm}

Now we consider set $DS_{a,b} = \{\alpha_0, \ldots, \alpha_{n-b}\}$. Clearly this set is closed under the addition. Let $\alpha_k, \alpha_l \in DS_{a,b}$ and $k, l \neq n - b$. Then from Proposition 3.2 follows that $\alpha_k \cdot \alpha_l = \alpha_0$. Also we have $\alpha_{n-b}\cdot \alpha_s = \alpha_0$, $\alpha_k \cdot \alpha_{n-b} = \alpha_k$ and $\alpha_{n-b} \cdot \alpha_{n-b} = \alpha_{n-b}$. Thus we prove

\vspace{3mm}

\textbf{Proposition 4.1 } \textsl{For any $a, b \in \mathcal{C}_n$  set $DS_{a,b}$ is a subsemiring of  semiring $S_{a,b}$.}

\vspace{3mm}

Note that in  semiring $DS_{a,b}$  endomorphism $\alpha_0$ is the zero element and  endomorphism $\alpha_{n-b}$ is the unique right identity.
Semiring $DS_{a,b}$ consist of all a -- nilpotent endomorphisms and the least idempotent endomorphism.

\vspace{3mm}

Now we define a mapping $D : \mathcal{STR}\{a,b\} \rightarrow \mathcal{STR}\{a,b\}$ by the rules $D(\alpha_k) = \alpha_{k-1}$ for any $k = 1, \ldots, n$, and $D(\alpha_0) = \alpha_0$.

Let $\alpha_k, \alpha_{\ell} \in \mathcal{STR}\{a,b\}$ and $k > \ell$. Then $D(\alpha_k + \alpha_{\ell}) = D(\alpha_k) = \alpha_{k-1} = \alpha_{k-1} + \alpha_{\ell-1} = D(\alpha_k) + D(\alpha_{\ell})$, that means $D$ is a linear mapping.

Let $\alpha_k, \alpha_{\ell} \in N_a$. Then from Proposition 3.2. follows that $D(\alpha_k\alpha_{\ell}) = D(\alpha_0) = \alpha_0$ and also $D(\alpha_k)\alpha_{\ell} = \alpha_{k-1}\alpha_{\ell} = \alpha_0$, $\alpha_k D(\alpha_{\ell}) = \alpha_k\alpha_{\ell-1} = \alpha_0$. So, we have
$$D(\alpha_k\alpha_{\ell}) = D(\alpha_k)\alpha_{\ell} + \alpha_k D(\alpha_{\ell}).$$

Since the same equalities are hold if we replace $\alpha_k$ with $\alpha_{n-b}$, then it follows
$$D(\alpha_{n-b}\alpha_{\ell}) = D(\alpha_{n-b})\alpha_{\ell} + \alpha_{n-b} D(\alpha_{\ell}).$$

For any $\alpha_k \in N_a$ we compute $D(\alpha_k \alpha_{n-b}) = D(\alpha_k) = \alpha_{k-1}$, $D(\alpha_k) \alpha_{n-b}
= \alpha_{k-1}\alpha_{n-b} = \alpha_{k-1}$ and $\alpha_kD(\alpha_{n-b}) = \alpha_k\alpha_{n-b-1} = \alpha_0$. So, we have
$$D(\alpha_k\alpha_{n-b}) = D(\alpha_k)\alpha_{n-b} + \alpha_k D(\alpha_{n-b}).$$

Also we compute $D(\alpha^2_{n-b}) = D(\alpha_{n-b}) = \alpha_{n-b-1}$ and $D(\alpha_{n-b})\alpha_{n-b} = \alpha_{n-b-1}\alpha_{n-b} = \alpha_{n-b-1}$, ${\alpha_{n-b}D(\alpha_{n-b}) = \alpha_{n-b}\alpha_{n-b-1} = \alpha_0}$. So, it follows
$$D(\alpha^2_{n-b}) = D(\alpha_{n-b})\alpha_{n-b} + \alpha_{n-b} D(\alpha_{n-b}).$$

Thus, we prove

\vspace{3mm}

\textbf{Proposition 4.1} \textsl{For any $a, b \in \mathcal{C}_n$  mapping $D$ is a derivation in  semiring $DS_{a,b}$.}

\vspace{3mm}

Now we prove that there are not differential semirings (under  derivation $D$) containing $DS_{a,b}$, which are subsemirings of $\mathcal{STR}\{a,b\}$.

 Let $\alpha_k \in I\!d_{a,b}$, where $k > n -b$.
We compute $D(\alpha^2_{k}) = D(\alpha_{k}) = \alpha_{k-1}$ and from Proposition 3.2 follows $D(\alpha_{k})\alpha_{k} = \alpha_{k-1}\alpha_{k} = \alpha_{k-1}$ and $\alpha_{k}D(\alpha_{k}) = \alpha_{k}\alpha_{k-1} = \alpha_k$. This means that\break  $D(\alpha^2_{k}) \neq D(\alpha_{k})\alpha_{k} + \alpha_{k}D(\alpha_{k})$.

 Let $\alpha_k \in N_b$. We compute $D(\alpha^2_{k}) = D(\alpha_{n}) = \alpha_{n-1}$, $D(\alpha_{k})\alpha_{k} = \alpha_{k-1}\alpha_{k} = \alpha_{n}$ and $\ds \alpha_{k}D(\alpha_{k}) = \alpha_{k}\alpha_{k-1} = \left\{ \begin{array}{ll} \alpha_{n-a}, & \mbox{if}\; k = n - a\\ \alpha_n, & \mbox{if}\; k > n - a \end{array} \right.$. Hence $D(\alpha_{k})\alpha_{k} + \alpha_{k}D(\alpha_{k}) = \alpha_n \neq \alpha_{n-1} = D(\alpha^2_k)$.

 Thus we prove that $D$ is not a Jordan derivation, see [Herst], in any subset of $\mathcal{STR}\{a,b\}$, which contains $DS_{a,b}$, that means $D$ is not a derivation. So, follows

\vspace{3mm}

\textbf{Proposition 4.2} \textsl{The semiring $DS_{a,b}$ is the maximal differential subsemiring (under the derivation $D$) of  string $\mathcal{STR}\{a,b\}$.}

\vspace{3mm}

Note that each of  subsets $I_0 = \{\alpha_0\}$, $I_1 = \{\alpha_0,\alpha_1\}$, $\ldots$, $I_{n-b-1} = N_a$ is a subsemiring of $N_a$ with trivial multiplication $\alpha_k\alpha_{\ell} = \alpha_0$ for any $0 \leq k, \ell \leq n - b - 1$. Since these semirings $I_k$, $0 \leq k \leq n - b- 1$, are closed under derivation $D$, then $I_k$ are differential subsemirings of $DS_{a,b}$. But from Proposition 3.2. follows that $I_k$ are ideals in  semiring $DS_{a,b}$. Hence

\vspace{3mm}

\textbf{Proposition 4.3} \textsl{ In  semiring  $DS_{a,b}$  there is a chain of differential ideals}
$$I_0 \subset \cdots \subset I_k = \{\alpha_0, \ldots, \alpha_k\} \subset \cdots \subset I_{n-b-1} = N_a.$$

\vspace{3mm}

Let $R$ be an arbitrary differential semiring with derivation $d$ and $I$ is a differential ideal of $R$. We consider
$$\int_R I = \left\{x | x \in R, \exists n \in \mathbb{N}\cup\{0\}, d^n(x) \in I\right\}.$$

Let $\ds x, y \in \int_R I$ where $d^m(x) \in I$ and $d^n(y) \in I$. If $m < n$ we have $d^n(x) = d^{n-m}(d^m(x)) \in I$. Then, using that $d$ is a linear map, follows $d^n(x + y) = d^n(x) + d^n(y) \in I$, that is $\ds x + y \in \int_R I$.

On the other hand, $\ds d^{m+n}(xy) = \sum_{k=0}^{m + n} \binom{m + n}{k} d^{m+n -k}(x)d^k(y) \in I$, which means that $\ds xy \in \int_R I$.
It is clear that $\ds x \in \int_R I$ implies $\ds d(x) \in \int_R I$. Thus, we prove

\vspace{3mm}

\textbf{Proposition 4.4} \textsl{Let $R$ be a differential semiring with derivation $d$ and $I$ is a differential ideal of $R$. Then $\ds \int_R I$ is a differential subsemiring of $R$. }

\vspace{3mm}

Using Propositions 4.3 and 4.4 we can describe the "differential structure" of semiring $DS_{a,b}$.

\vspace{3mm}

\textbf{Corollary 4.5 } \textsl{{For every differential ideal $I_k$ of semiring $DS_{a,b}$, where $k = 0, \ldots, n - b - 1$, follows  $\ds DS_{a,b} = \int_{DS_{a,b}} I_k$}.}

\vspace{5mm}

\centerline{\large 5. \hspace{0.5mm}Derivations in  string $\mathcal{STR}\{a,b\}$}

\vspace{3mm}

Now we use the well known Jordan multiplication in associative rings to define some new derivations in  strings $\mathcal{STR}\{a,b\}$.

Let $\alpha \in \mathcal{STR}\{a,b\}$. We define a mapping $\delta_\alpha : \mathcal{STR}\{a,b\} \rightarrow \mathcal{STR}\{a,b\}$ by the rule
 $$\delta_\alpha (\alpha_k) = \alpha\alpha_{k} + \alpha_k\alpha \;\; \mbox{for any}\; k = 0, 1, \ldots, n.$$

The main result in this section is

\vspace{3mm}

\textbf{Theorem 5.1} \textsl{For any $a, b \in \mathcal{C}_n$ and arbitrary $\alpha \in \mathcal{STR}\{a,b\}$   mapping $\delta_\alpha$ is a derivation in string $\mathcal{STR}\{a,b\}$.}

\emph{Proof.}
From $\delta_\alpha (\alpha_k + \alpha_\ell) = \alpha (\alpha_k + \alpha_\ell) + (\alpha_k + \alpha_\ell)\alpha = \alpha\alpha_{k} + \alpha_k\alpha + \alpha\alpha_{\ell} + \alpha_\ell\alpha = \delta_\alpha(\alpha_k) + \delta_\alpha(\alpha_\ell)$ where $k, \ell \in \{0, \ldots, n\}$ follows that  mapping $\delta_\alpha$ is a linear.

Now we prove  equality
$$\delta_\alpha (\alpha_k\alpha_\ell) = \delta_\alpha(\alpha_k)\alpha_\ell + \alpha_k\delta_\alpha(\alpha_\ell). \eqno(1)$$

\emph{Case 1.} Let $\; \alpha \in N_a$.

$\bullet$ If $\alpha_k \in N_a$, then $\delta_\alpha(\alpha_k) = \alpha\alpha_k + \alpha_k\alpha = \alpha_0 + \alpha_0 = \alpha_0$.

$\bullet$ If $\alpha_k \in I\!d_{a,b}$, then $\delta_\alpha(\alpha_k) = \alpha\alpha_k + \alpha_k\alpha = \alpha + \alpha_0 = \alpha$.

$\bullet$ If $\alpha_k \in N_b$, then $\delta_\alpha(\alpha_k) = \alpha\alpha_k + \alpha_k\alpha = \alpha_n + \alpha_0 = \alpha_n$.

\vspace{1mm}

\textbf{1.1.} Let $\alpha_\ell \in N_a$. Then $\delta_\alpha (\alpha_k\alpha_\ell) = \delta_\alpha(\alpha_0) = \alpha_0$, $\delta_\alpha (\alpha_k)\alpha_\ell = \alpha_0$ and $\alpha_k\delta_\alpha(\alpha_\ell) = \alpha_k\alpha_0 = \alpha_0$. So, (1) holds.

\textbf{1.2.} Let $\alpha_\ell \in N_b$. Then $\delta_\alpha (\alpha_k\alpha_\ell) = \delta_\alpha(\alpha_n) = \alpha_n$, $\delta_\alpha (\alpha_k)\alpha_\ell = \alpha_n$ and $\alpha_k\delta_\alpha(\alpha_\ell) = \alpha_k\alpha_n = \alpha_n$. So, (1) holds.

\textbf{1.3.} Let $\alpha_k \in N_a, \alpha_\ell \in I\!d_{a,b}$. Then $\delta_\alpha (\alpha_k\alpha_\ell) = \delta_\alpha(\alpha_k) = \alpha_0$, $\delta_\alpha (\alpha_k)\alpha_\ell = \alpha_0\alpha_\ell = \alpha_0$ and $\alpha_k\delta_\alpha(\alpha_\ell) = \alpha_k\alpha_0 = \alpha_0$. So, (1) holds.

\textbf{1.4.} Let $\alpha_k, \alpha_\ell \in I\!d_{a,b}$. Then $\delta_\alpha (\alpha_k\alpha_\ell) = \delta_\alpha(\alpha_k) = \alpha$, $\delta_\alpha (\alpha_k)\alpha_\ell = \alpha\alpha_\ell = \alpha$ and $\alpha_k\delta_\alpha(\alpha_\ell) = \alpha_k\alpha = \alpha_0$. So, (1) holds.

\textbf{1.5.} Let $\alpha_k \in N_b, \alpha_\ell \in I\!d_{a,b}$. Then $\delta_\alpha (\alpha_k\alpha_\ell) = \delta_\alpha(\alpha_k) = \alpha_n$, $\delta_\alpha (\alpha_k)\alpha_\ell = \alpha_n\alpha_\ell = \alpha_n$ and $\alpha_k\delta_\alpha(\alpha_\ell) = \alpha_k\alpha = \alpha_0$. So, (1) holds.

\vspace{2mm}

\emph{Case 2.} Let $\; \alpha \in I\!d_{a,b}$.

$\bullet$ If $\alpha_k \in N_a$, then $\delta_\alpha(\alpha_k) = \alpha\alpha_k + \alpha_k\alpha = \alpha_0 + \alpha_k = \alpha_k$.

$\bullet$ If $\alpha_k \in I\!d_{a,b}$, then $\delta_\alpha(\alpha_k) = \alpha\alpha_k + \alpha_k\alpha = \alpha + \alpha_k$.

$\bullet$ If $\alpha_k \in N_b$, then $\delta_\alpha(\alpha_k) = \alpha\alpha_k + \alpha_k\alpha = \alpha_n + \alpha_k = \alpha_n$.

\vspace{1mm}

\textbf{2.1.} Let $\alpha_\ell \in N_a$. Then (1) holds after the same equalities like in 1.1.

\textbf{2.2.} Let $\alpha_\ell \in N_b$. Then (1) holds after the same equalities like in 1.2.

\textbf{2.3.} Let $\alpha_k \in N_a, \alpha_\ell \in I\!d_{a,b}$. Then $\delta_\alpha (\alpha_k\alpha_\ell) = \delta_\alpha(\alpha_k) = \alpha_k$, $\delta_\alpha (\alpha_k)\alpha_\ell = \alpha_k\alpha_\ell = \alpha_k$ and $\alpha_k\delta_\alpha(\alpha_\ell) = \alpha_k(\alpha + \alpha_\ell) = \alpha_k\alpha + \alpha_k\alpha_\ell = \alpha_k + \alpha_k = \alpha_k$. So, (1) holds.

\textbf{2.4.} Let $\alpha_k, \alpha_\ell \in I\!d_{a,b}$. Then $\delta_\alpha (\alpha_k\alpha_\ell) = \delta_\alpha(\alpha_k) = \alpha + \alpha_k$, $\delta_\alpha (\alpha_k)\alpha_\ell = (\alpha + \alpha_k)\alpha_\ell = \alpha + \alpha_k$ and $\alpha_k\delta_\alpha(\alpha_\ell) = \alpha_k(\alpha + \alpha_\ell)= \alpha_k$. So, (1) holds.

\textbf{2.5.} Let $\alpha_k \in N_b, \alpha_\ell \in I\!d_{a,b}$. Then $\delta_\alpha (\alpha_k\alpha_\ell) = \delta_\alpha(\alpha_k) = \alpha_n$, $\delta_\alpha (\alpha_k)\alpha_\ell = \alpha_n\alpha_\ell = \alpha_n$ and $\alpha_k\delta_\alpha(\alpha_\ell) = \alpha_k(\alpha + \alpha_l) = \alpha_k$. So, (1) holds.

\vspace{2mm}

\emph{Case 3.} Let $\; \alpha \in N_b$.

$\bullet$ If $\alpha_k \in N_a$, then $\delta_\alpha(\alpha_k) = \alpha\alpha_k + \alpha_k\alpha = \alpha_0 + \alpha_n = \alpha_n$.

$\bullet$ If $\alpha_k \in I\!d_{a,b}$, then $\delta_\alpha(\alpha_k) = \alpha\alpha_k + \alpha_k\alpha = \alpha + \alpha_n = \alpha_n$.

$\bullet$ If $\alpha_k \in N_b$, then $\delta_\alpha(\alpha_k) = \alpha\alpha_k + \alpha_k\alpha = \alpha_n + \alpha_n = \alpha_n$.

\vspace{1mm}

\textbf{3.1.} Let $\alpha_\ell \in N_a$. Then $\delta_\alpha (\alpha_k\alpha_\ell) = \delta_\alpha(\alpha_0) = \alpha_n$, $\delta_\alpha (\alpha_k)\alpha_\ell = \alpha_0$ and $\alpha_k\delta_\alpha(\alpha_\ell) = \alpha_k\alpha_n = \alpha_n$. So, (1) holds.

\textbf{3.2.} Let $\alpha_\ell \in N_b$. Then (1) holds after the same equalities like in 1.2.

\textbf{3.3.} Let $\alpha_k \in N_a, \alpha_\ell \in I\!d_{a,b}$. Then $\delta_\alpha (\alpha_k\alpha_\ell) = \delta_\alpha(\alpha_k) = \alpha_n$, $\delta_\alpha (\alpha_k)\alpha_\ell = \alpha_n\alpha_\ell = \alpha_n$ and $\alpha_k\delta_\alpha(\alpha_\ell) = \alpha_k\alpha_n = \alpha_n$. So, (1) holds.

\textbf{3.4.} Let $\alpha_k, \alpha_\ell \in I\!d_{a,b}$. Then $\delta_\alpha (\alpha_k\alpha_\ell) = \delta_\alpha(\alpha_k) = \alpha_n$, $\delta_\alpha (\alpha_k)\alpha_\ell = \alpha_n\alpha_\ell = \alpha_n$ and $\alpha_k\delta_\alpha(\alpha_\ell) = \alpha_k\alpha_n = \alpha_n$. So, (1) holds.

\textbf{3.5.} Let $\alpha_k \in N_b, \alpha_\ell \in I\!d_{a,b}$. Then $\delta_\alpha (\alpha_k\alpha_\ell) = \delta_\alpha(\alpha_k) = \alpha_n$, $\delta_\alpha (\alpha_k)\alpha_\ell = \alpha_n\alpha_\ell = \alpha_n$ and $\alpha_k\delta_\alpha(\alpha_\ell) = \alpha_k\alpha_n = \alpha_n$. So, (1) holds and this completes the proof.

\vspace{3mm}

\textbf{Proposition 5.2 } \textsl{All the derivations $\delta_\alpha$, where $\alpha \in \mathcal{STR}\{a,b\}$ commutes.}

\emph{Proof.} Let $\alpha, \beta \in N_a$. For arbitrary $\alpha_k \in N_a$ follows $\delta_\alpha(\alpha_k) = \alpha_0$ and $\delta_\beta(\alpha_k) = \alpha_0$ what implies $\delta_\alpha \cdot \delta_\beta = \delta_\beta \cdot \delta_\alpha$. For arbitrary $\alpha_k \in I\!d_{a,b}$ we have $\delta_\alpha(\alpha_k) = \alpha$ and $\delta_\beta(\alpha_k) = \beta$. Then $\delta_\beta(\delta_\alpha(\alpha_k)) = \delta_\beta(\alpha) = \alpha_0$ and  $\delta_\alpha(\delta_\beta(\alpha_k)) = \delta_\alpha(\beta) = \alpha_0$ what implies $\delta_\alpha \cdot \delta_\beta = \delta_\beta \cdot \delta_\alpha$.

Let $\alpha \in N_a$, $\beta \in I\!d_{a,b}$. For arbitrary $\alpha_k \in N_a$ we have $\delta_\alpha(\alpha_k) = \alpha_0$ and $\delta_\beta(\alpha_k) = \alpha_k$. So, $\delta_\alpha(\delta_\beta(\alpha_k)) = \alpha_0 = \delta_\beta(\delta_\alpha(\alpha_k))$, i.e. $\delta_\alpha \cdot \delta_\beta = \delta_\beta \cdot \delta_\alpha$.  For arbitrary $\alpha_k \in I\!d_{a,b}$ we have $\delta_\alpha(\alpha_k) = \alpha_k$ and $\delta_\beta(\alpha_k) = \beta + \alpha_k$. Then $\delta_\beta(\delta_\alpha(\alpha_k)) = \beta + \alpha_k$ and $\delta_\alpha(\delta_\beta(\alpha_k)) = \delta_\alpha(\beta + \alpha_k) = \beta + \alpha_k$, that is $\delta_\alpha \cdot \delta_\beta = \delta_\beta \cdot \delta_\alpha$.

Let $\alpha, \beta \in I\!d_{a,b}$. For arbitrary $\alpha_k \in N_a$ follows $\delta_\alpha(\alpha_k) = \alpha_k$ and $\delta_\beta(\alpha_k) = \alpha_k$ what implies $\delta_\alpha \cdot \delta_\beta = \delta_\beta \cdot \delta_\alpha$.  For arbitrary $\alpha_k \in I\!d_{a,b}$ we have $\delta_\alpha(\alpha_k) = \alpha + \alpha_k$ and $\delta_\beta(\alpha_k) = \beta + \alpha_k$. Then $\delta_\beta(\delta_\alpha(\alpha_k)) = \delta_\beta(\alpha + \alpha_k) = \beta + \alpha + \alpha_k$ and $\delta_\alpha(\delta_\beta(\alpha_k)) = \delta_\alpha(\beta + \alpha_k) = \alpha + \beta + \alpha_k$, that is $\delta_\alpha \cdot \delta_\beta = \delta_\beta \cdot \delta_\alpha$.

Let  $\beta \in N_b$. For arbitrary $\alpha_k$ follows $\delta_\beta(\alpha_k) = \alpha_n$. Then for any $\alpha \in \mathcal{STR}\{a,b\}$ we have $\delta_\alpha(\delta_\beta(\alpha_k)) = \delta_\alpha(\alpha_n) = \alpha_n$ and $\delta_\beta(\delta_\alpha(\alpha_k) = \alpha_n$, i.e. $\delta_\alpha \cdot \delta_\beta = \delta_\beta \cdot \delta_\alpha$.

When $\alpha_k \in N_b$, then $\delta_\alpha(\alpha_k) = \alpha_n$ and $\delta_\beta(\alpha_k) = \alpha_n$ for any $\alpha, \beta \in \mathcal{STR}\{a,b\}$ what implies $\delta_\alpha \cdot \delta_\beta = \delta_\beta \cdot \delta_\alpha$.

\vspace{3mm}

Now we ask what structure has the set $\Delta$ of all derivations $\delta\alpha$, where  $\alpha \in \mathcal{STR}\{a,b\}$. First, we consider two examples.

\vspace{3mm}

\textbf{Example 5.3 } Let us fix $a = 0$ and examine  string $\mathcal{STR}\{0,b\}$. Consider the set of derivations $\Delta = \{\delta_\alpha \; |\; \alpha \in  I\!d_{0,b}\}$. Here only $\alpha_n$ is an element of   semiring $N_b$. For any $b$, $0 < b \leq n - 1$, we  compute $\delta_{\alpha_{n-b}}(\alpha_k) = \alpha_k$ for all $\alpha_k \in N_a$. It is clear that $\delta_{\alpha_{n-b}}(\alpha_{n-b}) = \alpha_{n-b}$ and $\delta_{\alpha_{n-b}}(\alpha_k) = \alpha_{n-b}\alpha_k + \alpha_k\alpha_{n-b} = \alpha_{n-b} + \alpha_k = \alpha_k$ for all $\alpha_k \in I\!d_{0,b}$. Moreover, $\delta_{\alpha_{n-b}}(\alpha_n) = \alpha_n$. Hence, $\delta_{\alpha_{n-b}}$ is an identity map.

Let $ \ell, m \in \{n-b, \ldots, n-1\}$, that is $\alpha_\ell, \alpha_m \in I\!d_{0,b}$, and $\ell \leq m$. Then follows $\delta_{\alpha_{\ell}}(\alpha_k) = \delta_{\alpha_{m}}(\alpha_k) = \alpha_k$ for all $\alpha_k \in N_a$ and also $\delta_{\alpha_{\ell}}(\alpha_n) = \delta_{\alpha_{m}}(\alpha_n) = \alpha_n$. For arbitrary $\alpha_k \in I\!d_{0,b}$ we have $\delta_{\alpha_{\ell}}(\alpha_k) = \alpha_\ell$, where $\alpha_k \leq \alpha_\ell$ and $\delta_{\alpha_{\ell}}(\alpha_k) = \alpha_k$ for any $k \in \{\ell + 1, \ldots, n - 1\}$. Hence, $\delta_{\alpha_m}(\delta_{\alpha_\ell}(\alpha_k)) = \delta_{\alpha_m}(\alpha_k)$ for every $\alpha_k \in I\!d_{0,b}$. Thus, $\delta_{\alpha_m}\cdot \delta_{\alpha_\ell} = \delta_{\alpha_m}$ when $\ell \leq m$.

 So, we conclude that the set of derivations $\Delta = \{\delta_{\alpha_{n-b}}, \ldots, \delta_{\alpha_{n-1}}\}$ is a commutative idempotent  semigroup, i.e. semilattice with identity $\delta_{\alpha_{n-b}}$.

\vspace{3mm}

\textbf{Example 5.4}
Let us fix $a = n-2$ and $b = n-1$. Now we compute the values of derivations:

$$\delta_{\alpha_0}(\alpha_0) = \alpha_0,\; \delta_{\alpha_0}(\alpha_1) = \alpha_0,\; \delta_{\alpha_0}(\alpha_2) = \alpha_n,\; \ldots,\; \delta_{\alpha_0}(\alpha_n) = \alpha_n,$$
$$\delta_{\alpha_1}(\alpha_0) = \alpha_0,\; \delta_{\alpha_1}(\alpha_1) = \alpha_1,\; \delta_{\alpha_1}(\alpha_2) = \alpha_n,\; \ldots,\; \delta_{\alpha_1}(\alpha_n) = \alpha_n,$$
$$\delta_{\alpha_2}(\alpha_0) = \alpha_n,\; \delta_{\alpha_2}(\alpha_1) = \alpha_n,\; \delta_{\alpha_2}(\alpha_2) = \alpha_n,\; \ldots,\; \delta_{\alpha_2}(\alpha_n) = \alpha_n,$$
$$ ...............................................................................................$$
$$\delta_{\alpha_n}(\alpha_0) = \alpha_n,\; \delta_{\alpha_n}(\alpha_1) = \alpha_n,\; \delta_{\alpha_n}(\alpha_2) = \alpha_n,\; \ldots,\; \delta_{\alpha_n}(\alpha_n) = \alpha_n.$$

Hence, $\delta_{\alpha_2} = \cdots = \delta_{\alpha_n}$. It is easy to verify that $\delta_{\alpha_0}^2 = \delta_{\alpha_0}$, $\delta_{\alpha_0}\cdot \delta_{\alpha_1} = \delta_{\alpha_0}$,  $\delta_{\alpha_0}\cdot \delta_{\alpha_2} = \delta_{\alpha_2}$, $\delta_{\alpha_1}^2 = \delta_{\alpha_1}$, $\delta_{\alpha_1}\cdot \delta_{\alpha_2} = \delta_{\alpha_2}$ and $\delta_{\alpha_2}^2 = \delta_{\alpha_2}$.

  So, we conclude that the  set of all derivations in  string $\mathcal{STR}\{n-2,n-1\}$ is the commutative idempotent semigroup $\Delta = \{\delta_{\alpha_{0}}, \delta_{\alpha_{1}}, \delta_{\alpha_{2}}\}$ with  multiplication table
$$\begin{array}{c|ccc} \cdot & \delta_{\alpha_{0}} & \delta_{\alpha_{1}} & \delta_{\alpha_{2}} \\ \hline
\delta_{\alpha_{0}} & \delta_{\alpha_{0}} & \delta_{\alpha_{0}} & \delta_{\alpha_{2}}\\
\delta_{\alpha_{1}} & \delta_{\alpha_{0}} & \delta_{\alpha_{1}} & \delta_{\alpha_{2}}\\
\delta_{\alpha_{2}} & \delta_{\alpha_{2}} & \delta_{\alpha_{2}} & \delta_{\alpha_{2}}
\end{array}.$$

\vspace{3mm}

\textbf{Theorem 5.5 } \textsl{For any $a, b \in \mathcal{C}_n$ the set of derivations $\Delta = \{\,\delta_{\alpha_0}, \ldots, \delta_{\alpha_{n-a}}\,\}$ in  string $\mathcal{STR}\{a,b\}$ is a semilattice with an identity $\delta_{n-b}$ and an  absorbing element $\delta_{n-a}$.
}

\emph{Proof.} Using Proposition 3.2 and  reasonings similar to those in proof of Theorem 5.1 we consider three cases.

\emph{Case 1.} Let $\; \alpha_\ell \in N_a$. It follows

$\bullet$ $\alpha_\ell \alpha_k = \alpha_k \alpha_\ell = \alpha_0$, where $\alpha_k \in N_a$ and then $\delta_{\alpha_\ell}(\alpha_k) = \alpha_0$;

$\bullet$ $\alpha_\ell \alpha_k = \alpha_\ell$ and $\alpha_k \alpha_\ell = \alpha_0$, where $\alpha_k \in I\!d_{a,b}$ and then $\delta_{\alpha_\ell}(\alpha_k) = \alpha_\ell$;

$\bullet$ $\alpha_\ell \alpha_k = \alpha_n$ and $\alpha_k \alpha_\ell = \alpha_0$, where $\alpha_k \in N_{b}$ and then $\delta_{\alpha_\ell}(\alpha_k) = \alpha_n$.

\emph{Case 2.} Let $\; \alpha_\ell \in I\!d_{a,b}$. It follows

$\bullet$ $\alpha_\ell \alpha_k = \alpha_0$ and $\alpha_k \alpha_\ell = \alpha_k$, where $\alpha_k \in N_a$ and then $\delta_{\alpha_\ell}(\alpha_k) = \alpha_k$;

$\bullet$ $\alpha_\ell \alpha_k = \alpha_\ell$ and $\alpha_k \alpha_\ell = \alpha_k$, where $\alpha_k \in I\!d_{a,b}$ and then $\delta_{\alpha_\ell}(\alpha_k) = \alpha_\ell + \alpha_k $;

$\bullet$ $\alpha_\ell \alpha_k = \alpha_n$ and $\alpha_k \alpha_\ell = \alpha_k$, where $\alpha_k \in N_{b}$ and then $\delta_{\alpha_\ell}(\alpha_k) = \alpha_n$.

\emph{Case 3.} Let $\; \alpha_\ell \in N_b$. It follows

$\bullet$ $\alpha_\ell \alpha_k = \alpha_0$ and $\alpha_k \alpha_\ell = \alpha_n$, where $\alpha_k \in N_a$ and then $\delta_{\alpha_\ell}(\alpha_k) = \alpha_n$;

$\bullet$ $\alpha_\ell \alpha_k = \alpha_\ell$ and $\alpha_k \alpha_\ell = \alpha_n$, where $\alpha_k \in I\!d_{a,b}$ and then $\delta_{\alpha_\ell}(\alpha_k) =  \alpha_n $;

$\bullet$ $\alpha_\ell \alpha_k = \alpha_k \alpha_\ell = \alpha_n$, where $\alpha_k \in N_{b}$ and then $\delta_{\alpha_\ell}(\alpha_k) = \alpha_n$.

From the last case we can conclude that $\delta_{n-a} = \cdots = \delta_n$.

Using the equalities in the cases above follows
$$\delta_{\alpha_0}(\alpha_0) = \alpha_0, \ldots, \delta_{\alpha_0}(\alpha_{n-b-1}) = \alpha_0,\delta_{\alpha_0}(\alpha_{n-b}) = \alpha_0, \ldots, \delta_{\alpha_0}(\alpha_{n-a-1}) = \alpha_0,$$ $$\delta_{\alpha_0}(\alpha_{n-a}) = \alpha_n, \ldots, \delta_{\alpha_0}(\alpha_n) = \alpha_n,$$
$$\delta_{\alpha_1}(\alpha_0) = \alpha_0, \ldots, \delta_{\alpha_1}(\alpha_{n-b-1}) = \alpha_0,\delta_{\alpha_1}(\alpha_{n-b}) = \alpha_1, \ldots, \delta_{\alpha_1}(\alpha_{n-a-1}) = \alpha_1,$$ $$\delta_{\alpha_1}(\alpha_{n-a}) = \alpha_n, \ldots, \delta_{\alpha_1}(\alpha_n) = \alpha_n,$$
$$............................................................................................................................................$$
$$\delta_{\alpha_{n-b-1}}(\alpha_0) = \alpha_0, \ldots, \delta_{\alpha_{n-b-1}}(\alpha_{n-b-1}) = \alpha_0,\delta_{\alpha_{n-b-1}}(\alpha_{n-b}) = \alpha_{n-b-1}, \ldots,$$ $$\delta_{\alpha_{n-b-1}}(\alpha_{n-a-1}) = \alpha_{n-b-1}, \delta_{\alpha_{n-b-1}}(\alpha_{n-a}) = \alpha_n, \ldots, \delta_{\alpha_{n-b-1}}(\alpha_n) = \alpha_n,$$
$$\delta_{\alpha_{n-b}}(\alpha_0) = \alpha_0, \delta_{\alpha_{n-b}}(\alpha_1) = \alpha_1 \ldots, \delta_{\alpha_{n-b}}(\alpha_{n-b-1}) = \alpha_{n-b-1},\delta_{\alpha_{n-b}}(\alpha_{n-b}) = \alpha_{n-b}, \ldots,$$ $$\delta_{\alpha_{n-b}}(\alpha_{n-a-1}) = \alpha_{n-a-1}, \delta_{\alpha_{n-b}}(\alpha_{n-a}) = \alpha_n, \ldots, \delta_{\alpha_{n-b}}(\alpha_n) = \alpha_n,$$
$$............................................................................................................................................$$
$$\delta_{\alpha_{n-a-1}}(\alpha_0) = \alpha_0, \delta_{\alpha_{n-a-1}}(\alpha_1) = \alpha_1, \ldots, \delta_{\alpha_{n-a-1}}(\alpha_{n-b-1}) = \alpha_{n-b-1},\delta_{\alpha_{n-a-1}}(\alpha_{n-b}) = \alpha_{n-a-1}, \ldots,$$ $$\delta_{\alpha_{n-a-1}}(\alpha_{n-a-1}) = \alpha_{n-a-1}, \delta_{\alpha_{n-a-1}}(\alpha_{n-a}) = \alpha_n, \ldots, \delta_{\alpha_{n-a-1}}(\alpha_n) = \alpha_n,$$
$$\delta_{\alpha_{n-a}}(\alpha_0) = \alpha_n, \ldots, \delta_{\alpha_{n-a}}(\alpha_{n-b-1}) = \alpha_n,\delta_{\alpha_{n-a}}(\alpha_{n-b}) = \alpha_n, \ldots, \delta_{\alpha_{n-a}}(\alpha_{n-a-1}) = \alpha_n,$$ $$\delta_{\alpha_{n-a}}(\alpha_{n-a}) = \alpha_n, \ldots, \delta_{\alpha_{n-a}}(\alpha_n) = \alpha_n.$$

Let us consider $\delta_{\alpha_\ell}$, where $\ell \leq n - a - 1$. Then for $k \leq n-a-1$ we compute $\delta_{\alpha_\ell}(\delta_{\alpha_0}(\alpha_k)) = \delta_{\alpha_\ell}(\alpha_0) = \alpha_0$ and for $k \geq n -a$ follows $\delta_{\alpha_\ell}(\delta_0(\alpha_k)) = \delta_{\alpha_\ell}(\alpha_n) = \alpha_n$. Also we compute $\delta_{\alpha_{n-2}}(\delta_{\alpha_0}(\alpha_k)) = \alpha_n$ for arbitrary $\alpha_k$.
Thus we prove that $\delta_{\alpha_0}\delta_{\alpha_\ell} = \delta_{\alpha_0}$ and $\delta_{\alpha_0}\delta_{\alpha_{n-a}} = \delta_{\alpha_{n-a}}$.

We find $\delta_{\alpha_1}(\delta_{\alpha_1}(\alpha_k)) = \delta_{\alpha_1}(\alpha_1) = \alpha_0$ where $k \leq n - a -1$ and $\delta_{\alpha_1}(\delta_{\alpha_1}(\alpha_k)) = \delta_{\alpha_1}(\alpha_n) = \alpha_n$ for $n - a \leq k \leq n$. So, we prove $\delta^2_{\alpha_1} = \delta_{\alpha_0}$.

Now we compute for $2 \leq \ell \leq n - b - 1$   elements $\delta_{\alpha_\ell}(\delta_{\alpha_1}(\alpha_0)) = \delta_{\alpha_\ell}(\alpha_0) = \alpha_0$
and $\delta_{\alpha_\ell}(\delta_{\alpha_1}(\alpha_1)) = \delta_{\alpha_\ell}(\alpha_0) = \alpha_0$, and also for $2 \leq k \leq n- a -1$  elements
$\delta_{\alpha_\ell}(\delta_{\alpha_1}(\alpha_k)) = \delta_{\alpha_\ell}(\alpha_1) = \alpha_0$. We have
$\delta_{\alpha_\ell}(\delta_{\alpha_1}(\alpha_{k}) = \delta_{\alpha_\ell}(\alpha_n) = \alpha_n$ for all $k \geq n-a$.
So, we prove $\delta_{\alpha_1}\delta_{\alpha_\ell} = \delta_{\alpha_0}$ for all $2 \leq \ell \leq n - b - 1$.

Now we compute for $n - b \leq \ell \leq n - a - 1$   elements $\delta_{\alpha_\ell}(\delta_{\alpha_1}(\alpha_0)) = \delta_{\alpha_\ell}(\alpha_0) = \alpha_0$ and $\delta_{\alpha_\ell}(\delta_{\alpha_1}(\alpha_1)) = \delta_{\alpha_\ell}(\alpha_0) = \alpha_0$. Also for $2 \leq k \leq n- a -1$  elements $\delta_{\alpha_\ell}(\delta_{\alpha_1}(\alpha_k)) = \delta_{\alpha_\ell}(\alpha_1) = \alpha_1$. Similarly, we have $\delta_{\alpha_\ell}(\delta_{\alpha_1}(\alpha_{k}) = \delta_{\alpha_\ell}(\alpha_n) = \alpha_n$ for all $k \geq n-a$. So, we prove $\delta_{\alpha_1}\delta_{\alpha_\ell} = \delta_{\alpha_1}$ for all $n-b \leq \ell \leq n - a - 1$.

For arbitrary $\alpha_k$ we compute $\delta_{\alpha_{n-a}}(\delta_{\alpha_1}(\alpha_{k})) =  \alpha_n$, so, we prove $\delta_{\alpha_1}\delta_{\alpha_{n-2}} = \delta_{\alpha_{n-2}}$.

Thus, using the similar and clear reasonings, and Proposition 5.2, we can construct the following table

$$\begin{array}{l|lllllllll}
    \cdot & \delta_{\alpha_0} & \delta_{\alpha_1} & \cdots & \delta_{\alpha_{n-b-1}} & \delta_{\alpha_{n-b}} & \delta_{\alpha_{n-b+1}}  &\cdots & \delta_{\alpha_{n-a-1}} & \delta_{\alpha_{n-a}}\\ \hline
     \delta_{\alpha_0} & \delta_{\alpha_0} & \delta_{\alpha_0} & \cdots & \delta_{\alpha_{0}} & \delta_{\alpha_{0}} & \delta_{\alpha_{0}}  &\cdots & \delta_{\alpha_{0}} & \delta_{\alpha_{n-a}}\\
     \delta_{\alpha_1} & \delta_{\alpha_0} & \delta_{\alpha_0} & \cdots & \delta_{\alpha_{0}} & \delta_{\alpha_{1}} & \delta_{\alpha_{1}}  &\cdots & \delta_{\alpha_{1}} & \delta_{\alpha_{n-a}}\\
       \vdots & \vdots & \vdots & \cdots & \vdots & \vdots & \vdots  &\cdots & \vdots & \vdots \\
     \delta_{\alpha_{n-b-1}} & \delta_{\alpha_0} & \delta_{\alpha_0} & \cdots & \delta_{\alpha_{0}} & \delta_{\alpha_{n-b-1}} & \delta_{\alpha_{n-b-1}}  &\cdots & \delta_{\alpha_{n-b-1}} & \delta_{\alpha_{n-a}}\\
     \delta_{\alpha_{n-b}} & \delta_{\alpha_0} & \delta_{\alpha_1} & \cdots & \delta_{\alpha_{n-b-1}} & \delta_{\alpha_{n-b}} & \delta_{\alpha_{n-b+1}}  &\cdots & \delta_{\alpha_{n-a-1}} & \delta_{\alpha_{n-a}}\\
    \delta_{\alpha_{n-b+1}} & \delta_{\alpha_0} & \delta_{\alpha_1} & \cdots & \delta_{\alpha_{n-b-1}} & \delta_{\alpha_{n-b+1}} & \delta_{\alpha_{n-b+1}}  &\cdots & \delta_{\alpha_{n-a-1}} & \delta_{\alpha_{n-a}}\\
     \vdots & \vdots & \vdots & \cdots & \vdots & \vdots & \vdots  &\cdots & \vdots & \vdots \\
    \delta_{\alpha_{n-b-1}} & \delta_{\alpha_0} & \delta_{\alpha_1} & \cdots & \delta_{\alpha_{n-b-1}} & \delta_{\alpha_{n-a-1}} & \delta_{\alpha_{n-a-1}}  &\cdots & \delta_{\alpha_{n-a-1}} & \delta_{\alpha_{n-a}}\\
     \delta_{\alpha_{n-a}} & \delta_{\alpha_{n-a}} & \delta_{\alpha_{n-a}} & \cdots & \delta_{\alpha_{n-a}} & \delta_{\alpha_{n-a}} & \delta_{\alpha_{n-a}}  &\cdots & \delta_{\alpha_{n-a}} & \delta_{\alpha_{n-a}}
  \end{array}
$$

\vspace{2mm}

This completes the proof that  set $\Delta = \{\,\delta_{\alpha_0}, \ldots, \delta_{\alpha_{n-a}}\,\}$ is a semilattice with identity $\delta_{n-b}$  absorbing element $\delta_{n-a}$.

\vspace{3mm}

Let $R$ be a differential semiring with set of derivations $\Delta = \{\delta_1, \ldots, \delta_m\}$ and $I$ be a differential ideal of $R$ that is closed under each derivation $\delta_i \in \Delta$. For any $i = 1, \ldots, m$ we denote
$$\int_R^{\delta_i} I = \left\{x | x \in R, \exists n \in \mathbb{N}\cup\{0\}, \delta_i^n(x) \in I\right\}.$$

From Proposition 3.2 follows that the set $I = \{\alpha_0,\alpha_{n}\}$ is an ideal in the string $\mathcal{STR}\{a,b\}$. From the proof of Theorem 5.1 we conclude that $I$ is closed under each derivation $\delta_{\alpha_i}$ where $\alpha_i \in \mathcal{STR}\{a,b\}$.
 An easy consequence of Proposition 4.4 and the proof of Theorem 5.5 is the following

\vspace{3mm}

\textbf{Corollary 5.6 } \textsl{For any $a, b \in \mathcal{C}_n$ subsemirings of the string $R = \mathcal{STR}\{a,b\}$ are}
$$\int_R^{\delta_{\alpha_0}}\! I = \cdots = \int_R^{\delta_{\alpha_{n-b-1}}}\!\! I = \int_R^{\delta_{\alpha_{n-a}}}\! I = R,\;\;\; \int_R^{\delta_{\alpha_{n-b}}}\! I = \cdots \int_R^{\delta_{\alpha_{n-a-1}}}\! I = I,$$
\textsl{where $\delta_{\alpha_i} \in \Delta = \{\,\delta_{\alpha_0}, \ldots, \delta_{\alpha_{n-a}}\,\}$ and $I = \{\alpha_0,\alpha_{n}\}$.}

\vspace{5mm}

 \centerline{\large 6. \hspace{0.5mm}Strings of arbitrary type}

\vspace{3mm}

For any $a_1, \ldots, a_m \in \mathcal{C}_n$, where $a_1 < a_2 < \ldots < a_m$, $m = 2, \ldots, n$ set
$$ \mathcal{STR}\{a_1, \ldots, a_m\} = \bigcup_{i = 1}^{m-1} \mathcal{STR}\{a_i,a_{i+1}\}
$$
is called a \textbf{\emph{string of type m}}.

\vspace{3mm}

Let $\alpha \in \mathcal{STR}\{a_i,a_{i+1}\}$ and $\beta \in \mathcal{STR}\{a_j,a_{j+1}\}$.
If $i = j$, then either $\alpha \leq \beta$, or $\beta \leq \alpha$. If $i < j$, then
 $$\alpha \leq \wr\,a_{i+1}, \ldots, a_{i+1} \,\wr \leq \wr\, a_j, \ldots, a_j \,\wr \leq \beta.$$

Hence, string $ \mathcal{STR}\{a_1, \ldots, a_m\}$ is a $(m-1)n + 1$ -- element chain with the least element $\wr\,a_1, \ldots, a_1\,\wr$ and the biggest element $\wr\,a_m, \ldots, a_m\,\wr$, so this string  is closed under the addition of  semiring $\widehat{\mathcal{E}}_{\mathcal{C}_n}$.

On the other hand, for $\alpha \in \mathcal{STR}\{a_i,a_{i+1}\}$ and $\beta \in \mathcal{STR}\{a_j,a_{j+1}\}$, where $i < j$, it is easy to show that
$\alpha \cdot \beta  \in \mathcal{STR}\{a_j,a_{j+1}\}$. Thus we prove

\vspace{3mm}

\textbf{Proposition 6.1 } \textsl{For any $a_1, \ldots, a_m \in \mathcal{C}_n$  string $\mathcal{STR}\{a_1, \ldots, a_m\}$ is a subsemiring of  semiring $\widehat{\mathcal{E}}_{{\mathcal{C}_n}}$.}

\vspace{3mm}

Immediately follows

\vspace{3mm}

\textbf{Corollary 6.2 } \textsl{For arbitrary subset $\{b_1, \ldots, b_\ell\} \subseteq \{a_1, \ldots, a_m\}$  string $\mathcal{STR}\{b_1, \ldots, b_\ell\}$ is a subsemiring of  string $\mathcal{STR}\{a_1, \ldots, a_m\}$.}

\vspace{3mm}

The elements of  semiring $\mathcal{STR}\{a_1, \ldots, a_m\}$, using the notations from section 3, are $\alpha_k^{a_\ell,a_{\ell+1}}$, where $\ell = 1, \ldots, m-1$ and $k = 0,\ldots,n$ is the number of elements of $\mathcal{C}_n$ with image equal to $a_{\ell+1}$. We can simplify this notations if we replace $\alpha_k^{a_\ell,a_{\ell+1}}$ with $\alpha_{k,\ell}$, where $\ell = 1, \ldots, m-1$. This means that $\alpha_{k,\ell}$ is the element of $\mathcal{STR}\{a_\ell, a_{\ell + 1}\}$ defined in the same way as in section 3. But using these notations we must remember that
$$\alpha_{0,2} = \alpha_{n,1},  \cdots, \alpha_{0,\ell+1} = \alpha_{n,\ell}, \cdots, \alpha_{o,m} = \alpha_{n,m-1}.$$

\vspace{3mm}

The next proposition is a generalization of Proposition 3.2.

\vspace{3mm}

\textbf{Proposition 6.3 } \textsl{Let  $a_1, \ldots, a_m \in \mathcal{C}_n$ and}
$$\mathcal{STR}\{a_1, \ldots, a_m\} = \{\alpha_{0,1}, \ldots, \alpha_{n,1},\alpha_{1,2}, \ldots, \alpha_{n,2}, \cdots, \alpha_{n,m-1},\alpha_{1,m}, \cdots, \alpha_{n,m}\}.$$

\textsl{For $k = 0, \ldots, n$ for the endomorphisms $\alpha_{k,\ell} \in \mathcal{STR}\{a_\ell,a_{\ell+1}\}$ and $\alpha_{s,r} \in \mathcal{STR}\{a_r,a_{r+1}\}$ follows}
$$\begin{array}{lll}
\alpha_{k,\ell} \cdot \alpha_{s,r} = \alpha_{0,r}, & \mbox{\textsl{where}} & 0 \leq s \leq n - a_{\ell+1} - 1\\
\alpha_{k,\ell} \cdot \alpha_{s,r} = \alpha_{k,r}, & \mbox{\textsl{where}} & n-a_{\ell+1} \leq s \leq n - a_\ell - 1\\
\alpha_{k,\ell} \cdot \alpha_{s,r} = \alpha_{n,r}, & \mbox{\textsl{where}} & n-a_\ell \leq s \leq n
\end{array}.$$

\emph{Proof.} Let  $0 \leq s \leq n - a_{\ell+1} - 1$. For arbitrary $i \in \mathcal{C}_n$ follows
$$(\alpha_{k,\ell}, \cdot \alpha_{s,r})(i) = \alpha_{s,r}(\alpha_{k,\ell}(i)) = \left\{ \begin{array}{ll} \alpha_{s,r}(a_\ell),& \mbox{if}\; i \leq n -k-1\\ \alpha_{s,r}(a_{\ell+1}), & \mbox{if} \; i \geq n - k \end{array}\right.$$ $$= \left\{ \begin{array}{ll} a_r,& \mbox{if}\; i \leq n -k-1\\ a_r, & \mbox{if} \; i \geq n - k \end{array} = \alpha_{0,r}(i) \right. . $$
So, $\alpha_{k,\ell} \cdot \alpha_{s,r} = \alpha_{0,r}$.

\vspace{2mm}

Let  $n - a_{\ell+1} \leq s \leq n - a_\ell - 1$. For all $i \in \mathcal{C}_n$ follows
$$(\alpha_{k,\ell}, \cdot \alpha_{s,r})(i) = \alpha_{s,r}(\alpha_{k,\ell}(i)) = \left\{ \begin{array}{ll} \alpha_{s,r}(a_\ell),& \mbox{if}\; i \leq n -k-1\\ \alpha_{s,r}(a_{\ell+1}), & \mbox{if} \; i \geq n - k \end{array} =\right.$$ $$\left\{ \begin{array}{ll} a_r,& \mbox{if}\; i \leq n -k-1\\ a_{r + 1}, & \mbox{if} \; i \geq  n - k \end{array} = \alpha_{k,r}(i) \right. . $$
So, $\alpha_{k,\ell} \cdot \alpha_{s,r} = \alpha_{k,r}$.

\vspace{2mm}

Let  $n - a_{\ell} \leq s \leq n $. For all $i \in \mathcal{C}_n$ follows
$$(\alpha_{k\ell}, \cdot \alpha_{s,r})(i) = \alpha_{s,r}(\alpha_{k,\ell}(i)) = \left\{ \begin{array}{ll} \alpha_{s,r}(a_\ell),& \mbox{if}\; i \leq n -k-1\\ \alpha_{s,r}(a_{\ell+1}), & \mbox{if} \; i \geq n - k \end{array} =\right.$$ $$ \left\{ \begin{array}{ll} a_{r+1},& \mbox{if}\; i \leq n -k-1\\ a_{r + 1}, & \mbox{if} \; i \geq n - k \end{array} = \alpha_{n,r}(i) \right. . $$
So, $\alpha_{k,\ell} \cdot \alpha_{s,r} = \alpha_{n,r}$ and this completes the proof.

\vspace{3mm}

Endomorphisms $\alpha_{0,\ell}$ where $\ell = 1, \ldots, m$ and $\alpha_{n,m}$ are called constant endomorphisms in [8]. According to [13], we denote the set of all constant endomorphisms of $\mathcal{STR}\{a_1, \ldots, a_m\}$  by $\mathcal{CO(STR}\{a_1, \ldots, a_m\})$. From Proposition 3.1 of [13] follows

\vspace{3mm}

\textbf{Corollary 6.4 } \textsl{Set $\mathcal{CO(STR}\{a_1, \ldots, a_m\})$ is an ideal of  semiring  $\mathcal{STR}\{a_1, \ldots, a_m\}$.}

\vspace{6mm}

 \centerline{\large 7. \hspace{0.5mm}Derivations in strings of arbitrary type}

\vspace{3mm}

We now proceed with the construction of derivations in strings $\mathcal{STR}\{a_1, \ldots, a_m\}$.

Suppose that there was $\ell = 1, \ldots, m-1$ such as $a_{\ell+1} - a_\ell \geq 2$. For this $\ell$ the set of endomorphisms $\alpha_{s,\ell}$, where $n - a_{\ell+1} \leq s \leq n - a_\ell - 1$ has 2 or more elements. So, if we put $p = n - a_{\ell+1}$, then $\alpha_p$ and $\alpha_{p+1}$ are from the semiring of the idempotent endomorphism in  string $\mathcal{STR}\{a_\ell,a_{\ell+1}\}$. Now, using Proposition 6.3 we compute
$$ \alpha_{p+1,\ell} \cdot \alpha_{p,\ell} = \alpha_{p+1,\ell}, \; \alpha_{p,\ell} \cdot \alpha_{p+1,\ell} = \alpha_{p,\ell}.
$$

For  $r \geq \ell +1$ we compute
$$\alpha_{p,\ell} \cdot \alpha_{p,r} = \alpha_{p,r}, \; \alpha_{p,\ell} \cdot \alpha_{p+1,r} = \alpha_{p,r}, \; \alpha_{p+1,\ell} \cdot \alpha_{p,r} = \alpha_{p+1,r}.
$$

To find  composition $\alpha_{p,r} \cdot \alpha_{p,\ell}$ we consider two possibilities:

1. If $n - a_{r+1} \leq p \leq n - a_r - 1$, then follows $\alpha_{p,r} \cdot \alpha_{p,\ell} = \alpha_{p,\ell}$.

2. If $n - a_r \leq p \leq n$, then $\alpha_{p,r} \cdot \alpha_{p,\ell} = \alpha_{n,\ell}$.

For composition $\alpha_{p+1,r} \cdot \alpha_{p,\ell}$ there are two similar possibilities, so we have either $\alpha_{p+1,r} \cdot \alpha_{p,\ell} = \alpha_{p+1,\ell}$, or $\alpha_{p+1,r} \cdot \alpha_{p,\ell} = \alpha_{n,\ell}$.

\vspace{3mm}

Consider a mapping $\delta_\alpha : \mathcal{STR}{\{a_1,\ldots,a_m\}} \rightarrow \mathcal{STR}{\{a_1,\ldots,a_m\}}$ defined by the rule
 $$\delta_\alpha (\alpha_{k,\ell}) = \alpha\alpha_{k,\ell} + \alpha_{k,\ell}\alpha, \eqno{(2)}$$
where $\alpha$ and $\alpha_{k,\ell}$ are arbitrary elements of the string $\mathcal{STR}{\{a_1,\ldots,a_m\}}$.
\vspace{2mm}

Then we find $\delta_{\alpha_{p,\ell}}(\alpha_{p+1,\ell}) = \alpha_{p+1,\ell}$.

Now we compute $\delta_{\alpha_{p,\ell}}(\alpha_{p,r}) = \alpha_{p,\ell} \cdot \alpha_{p,r} + \alpha_{p,r} \cdot \alpha_{p,\ell} = \alpha_{p,r} + \alpha_{p,r} \cdot \alpha_{p,\ell}$. But $\alpha_{p,r} > \alpha_{n,\ell} \geq \alpha_{p,\ell}$, hence $\delta_{\alpha_{p,\ell}}(\alpha_{p,r}) = \alpha_{p,r}$.

In such a way $\delta_{\alpha_{p,\ell}}(\alpha_{p+1,r}) = \alpha_{p,\ell} \cdot \alpha_{p+1,r} + \alpha_{p+1,r} \cdot \alpha_{p,\ell} = \alpha_{p,r} + \alpha_{p+1,r} \cdot \alpha_{p,\ell}$. Using the inequalities $\alpha_{p,r} > \alpha_{n,\ell} \geq \alpha_{p+1,\ell}$ follows $\delta_{\alpha_{p,\ell}}(\alpha_{p+1,r}) = \alpha_{p,r}$.

So, we have $\delta_{\alpha_{p,\ell}}(\alpha_{p+1,\ell} \cdot \alpha_{p,r}) = \delta_{\alpha_{p,\ell}}(\alpha_{p+1,r}) = \alpha_{p,r}$.

On the other hand, follows $\delta_{\alpha_{p,\ell}}(\alpha_{p+1,\ell}) \cdot \alpha_{p,r} = \alpha_{p+1,\ell} \cdot \alpha_{p,r} = \alpha_{p+1,r}$ and also $\alpha_{p+1,\ell} \cdot \delta_{\alpha_{p,\ell}}(\alpha_{p,r}) = \alpha_{p+1,\ell} \cdot \alpha_{p,r} = \alpha_{p+1,r}$. Hence $$\delta_{\alpha_{p,\ell}}(\alpha_{p+1,\ell}) \cdot \alpha_{p,r} + \alpha_{p+1,\ell} \cdot \delta_{\alpha_{p,\ell}}(\alpha_{p,r}) = \alpha_{p+1,r} \neq \alpha_{p,r} = \delta_{\alpha_{p,\ell}}(\alpha_{p+1,\ell} \cdot \alpha_{p,r}).$$

Thus, we show that  mapping $\delta_\alpha$ defined by (2), where $\alpha \in \mathcal{STR}{\{a_1,\ldots,a_m\}}$, in the general case, is not a derivation.

Note that there is not  a counterexample when $r = \ell$ because $\delta_{\alpha_{p,\ell}}(\alpha_{p+1,\ell} \cdot \alpha_{p,\ell}) = \delta_{\alpha_{p,\ell}}(\alpha_{p+1,\ell}) = \alpha_{p+1,\ell} = \delta_{\alpha_{p,\ell}}(\alpha_{p+1,\ell}) \cdot \alpha_{p,\ell} + \alpha_{p+1,\ell} \cdot \delta_{\alpha_{p,\ell}}(\alpha_{p,\ell})$. So, the last arguments do not contradict  the proof of Theorem 5.1.

To avoid possibility $a_{\ell+1} - a_\ell \geq 2$ we fix $m = n$ and $a_1 = 0, a_2 = 1, \ldots, a_n = n-1$.  Now we denote $\mathcal{STR}{\{0, 1,\ldots, n-1}\} = \mathcal{STR}\left(\widehat{\mathcal{E}}_{\mathcal{C}_n}\right)$.

\vspace{3mm}

\textbf{Example 7.1 } In  string $\mathcal{STR}\left(\widehat{\mathcal{E}}_{\mathcal{C}_4}\right)$ we consider endomorphisms $\alpha_{0,2} = \wr\, 2,2,2,2\,\wr$, $\alpha_{2,1} = \wr\, 0,0,1,1\,\wr$ and $\alpha_{2,3} = \wr\, 2,2,3,3\,\wr$. Since $\alpha_{2,1}\cdot \alpha_{2,3} = \alpha_{0,2}$ it follows that $\delta_{\alpha_{0,2}}(\alpha_{2,1}\cdot \alpha_{2,3}) = \delta_{\alpha_{0,2}}(\alpha_{0,2}) = \alpha_{0,2}$.  Now we compute $\delta_{\alpha_{0,2}}(\alpha_{2,1}) = \wr\, 2,2,2,2\,\wr \cdot \wr\, 0,0,1,1\,\wr + \wr\, 0,0,1,1\,\wr \cdot \wr\, 2,2,2,2\,\wr = \alpha_{0,2}$. Then
$\delta_{\alpha_{0,2}}(\alpha_{2,1}) \cdot \alpha_{2,3} = \alpha_{0,2} \cdot \alpha_{2,3} = \alpha_{0,3}$. Hence
$$\delta_{\alpha_{0,2}}(\alpha_{2,1}\cdot \alpha_{2,3}) \neq \delta_{\alpha_{0,2}}(\alpha_{2,1}) \cdot \alpha_{2,3} + \alpha_{2,1} \cdot \delta_{\alpha_{0,2}}(\alpha_{2,3}).$$

Thus, we showed that, even in the simplest case, in the string of type $m$, where $m > 2$, the mappings defined by (2) in general are not derivations.

\vspace{3mm}

Now we consider the ideal of constant endomorphisms.
It is clear that $\mathcal{CO}\left(\mathcal{STR}\left(\widehat{\mathcal{E}}_{\mathcal{C}_n}\right)\right) = \mathcal{CO}\left(\widehat{\mathcal{E}}_{\mathcal{C}_n}\right)$.
So, $\mathcal{CO}\left(\widehat{\mathcal{E}}_{\mathcal{C}_n}\right) = \{\kappa_0, \ldots, \kappa_{n-1}\}$ where $\kappa_i = \wr\,i, \ldots, i\,\wr$ for $i = 0, \ldots, n-1$.

Since $\kappa_i \cdot \kappa_j = \kappa_j$ it follows that $\delta_{\kappa_i}(\kappa_j) = \left\{\begin{array}{ll} \kappa_i,& \mbox{if}\; i \geq j\\    \kappa_j, & \mbox{if} \; i < j \end{array}\right.$.

More generally, for arbitrary $\alpha_{s,r} \in \mathcal{STR}\{r-1,r\}$, where $r = 1, \ldots, n - 1$ follows that $$\kappa_i \cdot \alpha_{s,r} = \left\{\begin{array}{ll} \kappa_{r-1},& \mbox{if}\; i \leq n-s-1\\ \kappa_{r}, & \mbox{if} \; i \geq n - s \end{array}\right. . \eqno{(3)} $$ Obviously $\alpha_{s,r}\cdot \kappa_i = \kappa_i$.

\vspace{2mm}

 Hence, there are three cases:

\textbf{1.} If $i \geq r$, then $\delta_{\alpha_{s,r}}(\kappa_i) = \kappa_i$.

\textbf{2.} If $i \leq r-1$ and $i \leq n -s-1$, then $\delta_{\alpha_{s,r}}(\kappa_i) = \kappa_{r-1}$.

\textbf{3.} If $i \leq r-1$ and $i \geq n -s$, then $\delta_{\alpha_{s,r}}(\kappa_i) = \kappa_{r}$.

\vspace{3mm}

We can now establish a result concerning the semiring $\mathcal{CO}\left(\widehat{\mathcal{E}}_{\mathcal{C}_n}\right)$ using Corollary 6.4.

\vspace{3mm}

\textbf{Proposition 7.2 } \textsl{Let $\delta_{\alpha_{s,r}} : \mathcal{CO}\left(\widehat{\mathcal{E}}_{\mathcal{C}_n}\right) \rightarrow \mathcal{CO}\left(\widehat{\mathcal{E}}_{\mathcal{C}_n}\right)$, where $\alpha_{s,r} \in \mathcal{STR}\left(\widehat{\mathcal{E}}_{\mathcal{C}_n}\right)$, be a mapping defined by $\delta_{\alpha_{s,r}} (\kappa_i) = \alpha_{s,r}\kappa_i + \kappa_i\alpha_{s,r}$
for arbitrary $\kappa_i \in \mathcal{CO}\left(\widehat{\mathcal{E}}_{\mathcal{C}_n}\right)$. Then for any $\kappa_i, \kappa_j \in \mathcal{CO}\left(\widehat{\mathcal{E}}_{\mathcal{C}_n}\right)$ follows}
$$\delta_{\alpha_{s,r}} (\kappa_i\cdot \kappa_j) = \delta_{\alpha_{s,r}} (\kappa_i)\cdot \kappa_j + \kappa_i \cdot \delta_{\alpha_{s,r}} (\kappa_j). \eqno{(4)}$$

\emph{Proof.} From  equality $\kappa_i \cdot \kappa_j = \kappa_j$  follows that $\delta_{\alpha_{s,r}} (\kappa_i\cdot \kappa_j) = \delta_{\alpha_{s,r}} (\kappa_j)$, $\kappa_i \cdot \delta_{\alpha_{s,r}}(\kappa_j) = \delta_{\alpha_{s,r}} (\kappa_j)$ and $\delta_{\alpha_{s,r}} (\kappa_i) \cdot \kappa_j = \kappa_j$. So (4) is equivalent to  equality $\delta_{\alpha_{s,r}} (\kappa_j) = \delta_{\alpha_{s,r}} (\kappa_j) + \kappa_j$. If $j \geq r$ from \textbf{1.} follows $\delta_{\alpha_{s,r}}(\kappa_j) = \kappa_j$, so, (4) holds. If $j \leq r-1$ and $j \leq n -s-1$, then from \textbf{ 2.} follows $\delta_{\alpha_{s,r}}(\kappa_j) = \kappa_{r-1}$ and from $\kappa_{r-1} = \kappa_{r-1} + \kappa_{j}$  equality (4) holds. When  $j \leq r-1$ and $j \geq n -s$, from \textbf{3.} follows $\delta_{\alpha_{s,r}}(\kappa_i) = \kappa_{r}$ which implies $\kappa_{r} = \kappa_{r} + \kappa_{j}$, so, (4) holds.

\vspace{3mm}

Is there a semiring which contain semiring $\mathcal{CO}\left(\widehat{\mathcal{E}}_{\mathcal{C}_n}\right)$ and is invariant under all mappings $\delta_{\alpha_{s,r}}$, where $\alpha_{s,r} \in \mathcal{STR}\left(\widehat{\mathcal{E}}_{\mathcal{C}_n}\right)$, so that  equality (3) holds for all elements of this semiring?

\vspace{3mm}

Studying this question, we consider the set $S = \mathcal{CO}\left(\widehat{\mathcal{E}}_{\mathcal{C}_n}\right) \cup \mathcal{STR}\{n-2,n-1\}$.

\vspace{3mm}

\textbf{Proposition 7.3 } \textsl{Set $S$ is a subsemiring of semiring $\mathcal{STR}\left(\widehat{\mathcal{E}}_{\mathcal{C}_n}\right)$ and is invariant under all mappings $\delta_{\alpha_{k,\ell}}$, where $\alpha_{k,\ell} \in \mathcal{STR}\left(\widehat{\mathcal{E}}_{\mathcal{C}_n}\right)$.
}

\emph{Proof.} Let $x, y \in S$. If either $x, y \in \mathcal{CO}\left(\widehat{\mathcal{E}}_{\mathcal{C}_n}\right)$, or $x,y \in  \mathcal{STR}\{n-2,n-1\}$, then from Corollary 6.4 and Proposition 3.1 follows that $x + y$ and $x\cdot y$ are from the same semiring.

Let $x = \alpha_{s,n-1} \in \mathcal{STR}\{n-2,n-1\}$ and $y = \kappa_i \in \mathcal{CO}\left(\widehat{\mathcal{E}}_{\mathcal{C}_n}\right)$, where $s, i  = 0, \ldots, n-1$. Then $\alpha_{s,n-1} + \kappa_i = \alpha_{s,n-1}$ if $i \leq n-2$ and $\alpha_{s,n-1} + \kappa_i = \kappa_i$ if $i = n-1$. Since $\alpha_{s,n-1}\cdot \kappa_i = \kappa_i$ and  $\kappa_i \cdot \alpha_{s,n-1} = \left\{\begin{array}{ll} \kappa_{n-2},& \mbox{if}\; i \leq n-s-1\\ \kappa_{n-1}, & \mbox{if} \; i \geq n - s \end{array}\right.$ follows that $x\cdot y, y\cdot x \in S$. So, we prove that $S$ is a subsemiring of $\mathcal{STR}\left(\widehat{\mathcal{E}}_{\mathcal{C}_n}\right)$.

Let $\alpha_{k,\ell} \in \mathcal{STR}\left(\widehat{\mathcal{E}}_{\mathcal{C}_n}\right)$. For any $\kappa_i \in \mathcal{CO}\left(\widehat{\mathcal{E}}_{\mathcal{C}_n}\right)$ from \textbf{1.}, \textbf{2.} and \textbf{3.} just before Proposition 7.2 follows that $\delta_{\alpha_{k,\ell}}(\kappa_i) \in \mathcal{CO}\left(\widehat{\mathcal{E}}_{\mathcal{C}_n}\right)$. Let $\alpha_{s,n-1} \in \mathcal{STR}\{n-2,n-1\}$. From Proposition 6.3 follows $\alpha_{k,\ell}\cdot \alpha_{s,n-1} = \kappa_{n-2}$ if $0 \leq s \leq n - \ell - 2$, $\alpha_{k,\ell}\cdot \alpha_{s,n-1} = \alpha_{k,n-1}$ if $s = n - \ell - 1$ and $\alpha_{k,\ell}\cdot \alpha_{s,n-1} = \kappa_{n-1}$ if $n - \ell \leq s \leq n - 1$. Since $\alpha_{s,n-1}\cdot \alpha_{k,\ell}$ is equal to $\kappa_{\ell-1}$, $\alpha_{s,\ell}$ and $\kappa_\ell$ in similar cases we can conclude that when $\ell \leq n-2$ follows
$$\delta_{\alpha_{k,\ell}}(\alpha_{s,n-1}) = \left\{\begin{array}{ll} \kappa_{n-2},& \mbox{if}\; 0 \leq s \leq n - \ell - 2\\ \alpha_{k,n-1}, & \mbox{if} \; s = n - \ell - 1\\ \kappa_{n-1}, & \mbox{if} \; n - \ell \leq s \leq n - 1 \end{array}.\right. \eqno{(5)}$$
If $\ell = n-1$ obviously $\delta_{\alpha_{k,n-1}}(\alpha_{s,n-1}) \in S$. So, $S$ is invariant under arbitrary $\delta_{\alpha_{k,\ell}}$ and this completes the proof.

\vspace{3mm}

The following counterexample shows that, in general,  map $\delta_{\alpha_{k,\ell}}$ is not a derivation.

\vspace{3mm}

\textbf{Example 7.4} In  string $\mathcal{STR}\left(\widehat{\mathcal{E}}_{\mathcal{C}_4}\right)$ we consider  endomorphisms      $\alpha_{3,2} = \wr\, 1,2,2,2\,\wr$, $\kappa_1 = \wr\, 1,1,1,1\,\wr \in \mathcal{CO}\left(\widehat{\mathcal{E}}_{\mathcal{C}_4}\right)$ and $\alpha_{2,3} = \wr\, 2,2,3,3\,\wr \in \mathcal{STR}\{2,3\}$.

Since $\kappa_1\cdot \alpha_{2,3} = \kappa_2$  follows that $$\delta_{\alpha_{3,2}}(\kappa_1\cdot \alpha_{2,3}) = \delta_{\alpha_{3,2}}(\kappa_2) = \wr\, 1,2,2,2\,\wr\cdot \wr\, 2,2,2,2\,\wr + \wr\, 2,2,2,2\,\wr \cdot \wr\, 1,2,2,2\,\wr = \wr\, 2,2,2,2\,\wr = \kappa_2.$$  Now we compute $\delta_{\alpha_{3,2}}(\kappa_1) = \wr\, 1,2,2,2\,\wr \cdot \wr\, 1,1,1,1\,\wr + \wr\, 1,1,1,1\,\wr \cdot \wr\, 1,2,2,2\,\wr = \kappa_2$. Then
$\delta_{\alpha_{3,2}}(\kappa_1) \cdot \alpha_{2,3} = \wr\, 2,2,2,2\,\wr \cdot \wr\, 2,2,3,3\,\wr = \wr\, 3,3,3,3\,\wr = \kappa_3$. Hence
$$\delta_{\alpha_{3,2}}(\kappa_1\cdot \alpha_{2,3}) \neq \delta_{\alpha_{3,2}}(\kappa_1) \cdot \alpha_{2,3} + \kappa_1 \cdot \delta_{\alpha_{3,2}}(\alpha_{2,3}).$$

\vspace{3mm}

Now we present a class of maps of type $\delta_\alpha$ which are derivations in the whole semiring $\mathcal{STR}\left(\widehat{\mathcal{E}}_{\mathcal{C}_n}\right)$.  We need some preliminary lemmas.

\vspace{3mm}

Using that $\delta_\alpha(\beta) = \delta_\beta(\alpha)$ and formulas (5) we obtain

\vspace{3mm}

\textbf{Lemma 7.5} \textsl{For any $\ell = 1, \ldots, n-2$ and $k, s = 0, \ldots, n$ follows}
$$\delta_{\alpha_{s,n-1}}(\alpha_{k,\ell}) = \left\{\begin{array}{ll} \kappa_{n-2}, & \mbox{if}\;\, 0 \leq s \leq n - \ell -1\\ \alpha_{k,n-1}, & \mbox{if}\;\, s = n - \ell\\ \kappa_{n-1}, & \mbox{if}\;\, n - \ell + 1 \leq s \leq n \end{array} . \right.$$

\vspace{3mm}

Let us make one necessary observation.

\vspace{3mm}

\textbf{Lemma 7.6} \textsl{For any $q = 1, \ldots, n-2$, $k, p, s = 0, \ldots, n$ and arbitrary $\ell$ follows}
$$\delta_{\alpha_{s,n-1}}(\alpha_{k,\ell})\cdot \alpha_{p,q} + \alpha_{k,\ell}\cdot \delta_{\alpha_{s,n-1}}(\alpha_{p,q}) = \alpha_{k,\ell}\cdot \delta_{\alpha_{s,n-1}}(\alpha_{p,q}).$$

\emph{Proof.} Let $\ell \leq n-2$. Using Lemma 7.5 follows  $\delta_{\alpha_{s,n-1}}(\alpha_{k,\ell}) \leq \kappa_{n-1}$. Then  $\delta_{\alpha_{s,n-1}}(\alpha_{k,\ell})\cdot \alpha_{p,q} \leq \kappa_{n-1}\cdot \alpha_{p,q} = \kappa_q$. The last equality follows from (3).

Let $\ell = n-1$. Then
$$\delta_{\alpha_{s,n-1}}(\alpha_{k,n-1}) = \left\{\begin{array}{ll} \kappa_{n-2}, & \mbox{if}\;\, k = s = 0,\; k = 0\; \mbox{and}\; s = 1,\; k = 1 \; \mbox{and}\; s = 0\\ \alpha_{1,n-1}, & \mbox{if}\;\, k = s = 1\\ \kappa_{n-1}, & \mbox{if}\;\, k \geq 2 \; \mbox{or}\; s \geq 2 \end{array} . \right. \eqno{(6)}$$

 Analogously, we have $\delta_{\alpha_{s,n-1}}(\alpha_{k,n-1})\cdot \alpha_{p,q} \leq \kappa_{n-1}\cdot \alpha_{p,q} = \kappa_q$.

 \vspace{2mm}

On the other hand, from Lemma 7.5 we have $\delta_{\alpha_{s,n-1}}(\alpha_{p,q}) \geq \kappa_{n-2}$. Then for arbitrary $\ell$ follows that $\alpha_{k,\ell}\cdot \delta_{\alpha_{s,n-1}}(\alpha_{p,q}) \geq \alpha_{k,\ell}\cdot \kappa_{n-2} = \kappa_{n-2}$. So, we have
$$\delta_{\alpha_{s,n-1}}(\alpha_{k,\ell})\cdot \alpha_{p,q} \leq \kappa_q \leq \kappa_{n-2} \leq \alpha_{k,\ell}\cdot \delta_{\alpha_{s,n-1}}(\alpha_{p,q}). $$

\vspace{3mm}

Now we shall prove the main result in this section.

\vspace{3mm}

\textbf{Theorem 7.7}\textsl{ For any $s \geq 2$ the map $\delta_{\alpha_{s,n-1}} : \mathcal{STR}\left(\widehat{\mathcal{E}}_{\mathcal{C}_n}\right) \rightarrow \mathcal{STR}\left(\widehat{\mathcal{E}}_{\mathcal{C}_n}\right)$ defined by  equality
$$\delta_{\alpha_{s,n-1}}(\alpha) = \alpha_{s,n-1}\cdot \alpha + \alpha\cdot \alpha_{s,n-1},$$
where $\alpha, \alpha_{s,n-1} \in \mathcal{STR}\left(\widehat{\mathcal{E}}_{\mathcal{C}_n}\right)$ is a derivation of $\mathcal{STR}\left(\widehat{\mathcal{E}}_{\mathcal{C}_n}\right)$. }

\emph{Proof.} \textbf{A.} First, we consider  endomorphisms $\alpha_{k,\ell}$ and $\alpha_{p,q}$ so that $\ell, q = 1, \ldots, n-2$. For $k, p = 0, \ldots, n$ we verify that
$$\delta_{\alpha_{s,n-1}}(\alpha_{k,\ell}\cdot \alpha_{p,q}) = \delta_{\alpha_{s,n-1}}(\alpha_{k,\ell})\cdot \alpha_{p,q} + \alpha_{k,\ell}\cdot \delta_{\alpha_{s,n-1}}(\alpha_{p,q}). \eqno{(7)}$$

 From Proposition 6.3 follows
$$ \begin{array}{ll} \alpha_{k,\ell}\cdot \alpha_{p,q} = \kappa_{q-1}, & \mbox{if}\;\, 0 \leq p \leq n - \ell -1\\ \alpha_{k,\ell}\cdot \alpha_{p,q} = \alpha_{k,q}, & \mbox{if}\;\, p = n - \ell\\ \alpha_{k,\ell}\cdot \alpha_{p,q} = \kappa_{q}, & \mbox{if}\;\, n - \ell + 1 \leq p \leq n \end{array} . \eqno{(8)}$$

\emph{Case 1.} Let $0 \leq s \leq n - q -1$. Then, using Lemma 7.5, for arbitrary $p$ follows
 $$\delta_{\alpha_{s,n-1}}(\alpha_{k,\ell}\cdot \alpha_{p,q}) = \kappa_{n-2}.$$

Now from Lemma 7.5. and Lemma 7.6 we find $\delta_{\alpha_{s,n-1}}(\alpha_{k,\ell})\cdot \alpha_{p,q} + \alpha_{k,\ell}\cdot \delta_{\alpha_{s,n-1}}(\alpha_{p,q}) = \alpha_{k,\ell}\kappa_{n-2} = \kappa_{n-2}$. So, equality (7) holds.
\vspace{1mm}

\emph{Case 2.} Let $s = n - q$.

\emph{2.1} If $0 \leq p \leq n - \ell - 1$, then from (8) follows $\alpha_{k,\ell}\cdot \alpha_{p,q} = \kappa_{q-1}$. Hence
 $$\delta_{\alpha_{s,n-1}}(\alpha_{k,\ell}\cdot \alpha_{p,q}) = \delta_{\alpha_{s,n-1}}(\kappa_{q-1}) = \kappa_{n-2}.$$

Using Lemma 7.6 first and then Lemma 7.5 we have $\delta_{\alpha_{s,n-1}}(\alpha_{k,\ell})\cdot \alpha_{p,q} + \alpha_{k,\ell}\cdot \delta_{\alpha_{s,n-1}}(\alpha_{p,q}) = \alpha_{k,\ell}\alpha_{p,n-1} = \kappa_{n-2}$. So, again equality (7) holds.

\vspace{1mm}

\emph{2.2} If $ p = n - \ell$, then from (6) follows $\alpha_{k,\ell}\cdot \alpha_{p,q} = \alpha_{k,q}$. Thus follows
 $$\delta_{\alpha_{s,n-1}}(\alpha_{k,\ell}\cdot \alpha_{p,q}) = \delta_{\alpha_{s,n-1}}(\alpha_{k,q}) = \alpha_{k,n-1}.$$

Now Lemma 7.5 and Lemma 7.6 yields $\delta_{\alpha_{s,n-1}}(\alpha_{k,\ell})\cdot \alpha_{p,q} + \alpha_{k,\ell}\cdot \delta_{\alpha_{s,n-1}}(\alpha_{p,q}) = \alpha_{k,\ell}\alpha_{p,n-1} = \alpha_{k,n-1}$. So, equality (7) holds.

\vspace{1mm}

\emph{2.3} If $n - \ell + 1 \leq p \leq n$, then from (8) follows $\alpha_{k,\ell}\cdot \alpha_{p,q} = \kappa_q$. So, we have
 $$\delta_{\alpha_{s,n-1}}(\alpha_{k,\ell}\cdot \alpha_{p,q}) = \delta_{\alpha_{s,n-1}}(\kappa_q) = \kappa_{n-1}.$$

 Now Lemma 7.5 and Lemma 7.6 imply $\delta_{\alpha_{s,n-1}}(\alpha_{k,\ell})\cdot \alpha_{p,q} + \alpha_{k,\ell}\cdot \delta_{\alpha_{s,n-1}}(\alpha_{p,q}) = \alpha_{k,\ell}\kappa_{n-1} = \kappa_{n-1}$. So, equality (7) holds.

\vspace{1mm}

\emph{Case 3.} Let $n - q + 1 \leq s \leq n$. Then, from Lemma 7.5, for arbitrary $p$ follows
 $$\delta_{\alpha_{s,n-1}}(\alpha_{k,\ell}\cdot \alpha_{p,q}) = \kappa_{n-1}.$$

From Lemma 7.5. and Lemma 7.6 we obtain $\delta_{\alpha_{s,n-1}}(\alpha_{k,\ell})\cdot \alpha_{p,q} + \alpha_{k,\ell}\cdot \delta_{\alpha_{s,n-1}}(\alpha_{p,q}) = \alpha_{k,\ell}\kappa_{n-1} = \kappa_{n-1}$. So, finally  equality (7) holds.

\vspace{2mm}

\textbf{B.} The second possibility is when $\alpha_{k,n-1}, \alpha_{p,n-1} \in \mathcal{STR}\{n-2,n-1\}$. Now from Theorem 5.1. follows that $\delta_{\alpha_{s,n-1}}$ satisfies the Leibnitz's rule.

\vspace{2mm}

\textbf{C.} Another possibility is to calculate $\delta_{\alpha_{s,n-1}}(\alpha_{k,\ell}\cdot \alpha_{p,n-1})$, where $\ell \leq n -2$. Now  equalities (8) for $q = n-1$ are
$$ \begin{array}{ll} \alpha_{k,\ell}\cdot \alpha_{p,n-1} = \kappa_{n-2}, & \mbox{if}\;\, 0 \leq p \leq n - \ell -1\\ \alpha_{k,\ell}\cdot \alpha_{p,n-1} = \alpha_{k,n-1}, & \mbox{if}\;\, p = n - \ell\\ \alpha_{k,\ell}\cdot \alpha_{p,n-1} = \kappa_{n-1}, & \mbox{if}\;\, n - \ell + 1 \leq p \leq n \end{array} . \eqno{(9)}$$

Hence, using that $s \geq 2$, for arbitrary $p$ we obtain $\delta_{\alpha_{s,n-1}}(\alpha_{k,\ell}\cdot \alpha_{p,n-1}) \geq \delta_{\alpha_{s,n-1}}(\kappa_{n-2}) = \kappa_{n-1}$. The last equality follows from (6).

From (6) and (9) we have $\alpha_{k,\ell} \cdot \delta_{\alpha_{s,n-1}}(\alpha_{p,n-1}) = \alpha_{k,\ell}\cdot \kappa_{n-1} = \kappa_{n-1}$.
Thus we prove that
$$\delta_{\alpha_{s,n-1}}(\alpha_{k,\ell}\cdot \alpha_{p,n-1}) = \delta_{\alpha_{s,n-1}}(\alpha_{k,\ell})\cdot \alpha_{p,n-1} + \alpha_{k,\ell}\cdot \delta_{\alpha_{s,n-1}}(\alpha_{p,n-1}).$$

\vspace{2mm}

\textbf{D.} The last possibility is the same as in \textbf{C.}, but we shall prove  equality
$$\delta_{\alpha_{s,n-1}}(\alpha_{p,n-1}\cdot \alpha_{k,\ell}) = \delta_{\alpha_{s,n-1}}(\alpha_{p,n-1})\cdot \alpha_{k,\ell} + \alpha_{p,n-1}\cdot \delta_{\alpha_{s,n-1}}(\alpha_{k,\ell}), \eqno{(10)}$$
where $\ell \leq n-2$.

Now from  equalities (8) we obtain
$$ \begin{array}{ll} \alpha_{p,n-1}\cdot \alpha_{k,\ell} = \kappa_{\ell - 1}, & \mbox{if}\;\, k = 0\\ \alpha_{p,n-1}\cdot \alpha_{k,\ell} = \alpha_{p,\ell}, & \mbox{if}\;\, k = 1\\ \alpha_{p,n-1}\cdot \alpha_{k,\ell} = \kappa_{\ell}, & \mbox{if}\;\, 2 \leq k \leq n \end{array} . \eqno{(11)}$$

\emph{Case 1.} Let $0 \leq s \leq n - \ell -1$. Then, using Lemma 7.5, for arbitrary $k$ follows
 $$\delta_{\alpha_{s,n-1}}(\alpha_{p,n-1}\cdot \alpha_{k,\ell}) = \kappa_{n-2}.$$

Now  Lemma 7.5. and Lemma 7.6 imply $\delta_{\alpha_{s,n-1}}(\alpha_{p,n-1})\cdot \alpha_{k,\ell} + \alpha_{p,n-1}\cdot \delta_{\alpha_{s,n-1}}(\alpha_{k,\ell}) = \alpha_{p,n-1}\kappa_{n-2} = \kappa_{n-2}$. So, equality (10) holds.

\vspace{1mm}

\emph{Case 2.} Let $s = n - \ell$.

\emph{2.1} Let $k = 0$. Now from (11) follows $\alpha_{p,n-1}\cdot \alpha_{k,\ell} = \kappa_{\ell - 1}$. Therefore, from Lemma 7.5, follows
 $$\delta_{\alpha_{s,n-1}}(\alpha_{p,n-1}\cdot \alpha_{k,\ell})  = \delta_{\alpha_{s,n-1}}(\kappa_{\ell-1}) = \kappa_{n-2}.$$

Using Lemma 7.6 first and then Lemma 7.5 we have $\delta_{\alpha_{s,n-1}}(\alpha_{p,n-1})\cdot \alpha_{k,\ell} + \alpha_{p,n-1}\cdot \delta_{\alpha_{s,n-1}}(\alpha_{k,\ell}) = \alpha_{p,n-1}\kappa_{n-2} = \kappa_{n-2}$. So, equality (10) holds.

\vspace{1mm}

\emph{2.2} Let $k = 1$. Now from (11) follows $\alpha_{p,n-1}\cdot \alpha_{k,\ell} = \alpha_{p,\ell}$. Thus
 $$\delta_{\alpha_{s,n-1}}(\alpha_{p,n-1}\cdot \alpha_{k,\ell})  = \delta_{\alpha_{s,n-1}}(\alpha_{p,\ell}) = \alpha_{p,n-1}.$$

 Now from Lemma 7.5  and  Lemma 7.6 we obtain $\delta_{\alpha_{s,n-1}}(\alpha_{p,n-1})\cdot \alpha_{k,\ell} + \alpha_{p,n-1}\cdot \delta_{\alpha_{s,n-1}}(\alpha_{k,\ell}) = \alpha_{p,n-1}\alpha_{k,n-1}  = \alpha_{p,n-1}$. So, equality (10) holds.
\vspace{1mm}

\emph{2.3} Let $2 \leq k \leq n$. From (11) follows $\alpha_{p,n-1}\cdot \alpha_{k,\ell} = \kappa_\ell$. Hence
 $$\delta_{\alpha_{s,n-1}}(\alpha_{p,n-1}\cdot \alpha_{k,\ell})  = \delta_{\alpha_{s,n-1}}(\kappa_\ell) = \kappa_{n-1}.$$

 Now  Lemma 7.5  and  Lemma 7.6 imply $\delta_{\alpha_{s,n-1}}(\alpha_{p,n-1})\cdot \alpha_{k,\ell} + \alpha_{p,n-1}\cdot \delta_{\alpha_{s,n-1}}(\alpha_{k,\ell}) = \alpha_{p,n-1}\alpha_{k,n-1}  = \kappa_{n-1}$. The last equality follows from the condition $k \geq 2$. So, equality (10) holds.

\vspace{1mm}

\emph{Case 3.} Let $n - \ell + 1 \leq s \leq n$. Then, using Lemma 7.5, for arbitrary $k$ follows
 $$\delta_{\alpha_{s,n-1}}(\alpha_{p,n-1}\cdot \alpha_{k,\ell}) = \kappa_{n-1}.$$

  From Lemma 7.5  and  Lemma 7.6 we find $\delta_{\alpha_{s,n-1}}(\alpha_{p,n-1})\cdot \alpha_{k,\ell} + \alpha_{p,n-1}\cdot \delta_{\alpha_{s,n-1}}(\alpha_{k,\ell}) = \alpha_{p,n-1}\kappa_{n-1}  = \kappa_{n-1}$. So, equality (10) holds and this completes the proof.

\vspace{3mm}

After Theorem 5.5 we consider  set $\ds \int_R^{\delta_i} I$, where $R$ is a differential semiring with a set of derivations $\Delta = \{\delta_1, \ldots, \delta_m\}$, $\delta_i \in \Delta$ and $I$ is a differential ideal of $R$, closed under each derivation $\delta_i$. From Proposition 4.4. follows that $\ds \int_R^{\delta_i} I$ is a subsemiring of $R$. In Theorem 7.7 we prove that string $R = \mathcal{STR}\left(\widehat{\mathcal{E}}_{\mathcal{C}_n}\right)$ is a differential semiring with set of derivations $\Delta = \{ \delta_{\alpha_{s,n-1}} | \alpha_{s,n-1} \in \mathcal{STR}\{n-2,n-1\}, \; 2 \leq s \leq n\}$. From Proposition 7.2 we know that $I = \mathcal{CO}\left(\widehat{\mathcal{E}}_{\mathcal{C}_n}\right)$ is a differential ideal of $R$, closed under all derivations of $\Delta$.

 \vspace{1mm}

 Now from Lemma 7.5 and equalities (6) easily follows

\vspace{3mm}

\textbf{Corollary 7.8 } \textsl{If $\delta_{\alpha_{s,n-1}} \in \Delta$ then $\ds \int_R^{\delta_{\alpha_{s,n-1}}} I = R$. }

\vspace{8mm}

\centerline{\large References}

\vspace{4mm}

{\small

[1] K. I. Beidar, W. S. Martindale III and A. V.  Mikhalev,  Rings with generalized identities, Marcel Dekker, 1996.

[2]  O. Ganyushkin and V. Mazorchuk, Classical Finite Transformation Semigroups: An Introduction,
Springer-Verlag London Limited, 2009.

[3]  J. Golan, Semirings and Their Applications, Kluwer, Dordrecht, 1999.

[4]  I. N. Herstein, Jordan derivations of prime rings, Proc. Amer. Math. Soc. 8 (1957) 1104-1110.

[5]  I. N. Herstein, Lie and Jordan structures in simple, associative rings, Bull. Amer. Math. Soc. 67 (6) (1961)
517-531.

[6]  I. N. Herstein,  Noncommutative Rings, Carus Mathematical Monographs, 1968.

[7]  I. N. Herstein, On the Lie structure of an associative ring, J. Algebra 14 (1970) 561-571.

[8]  J. Je$\hat{\mbox{z}}$ek, T. Kepka and M. Mar\`{o}ti, The endomorphism semiring of a semilattice,
Semigroup Forum, 78 (2009), 21 -- 26.

[9]  E. R. Kolchin, Differential Algebra and Algebraic Groups, Academic Press, New York, London, 1973.

[10]  C. Monico, On finite congruence-simple semirings, J. Algebra 271 (2004)
846 -- 854.

[11]  J. F. Ritt, Differential Algebra, Amer. Math. Soc. Colloq. Publ., Vol 33, New York, 1950.

[12]  I. Trendafilov and  D. Vladeva, The endomorphism semiring of a finite chain, Proc.
Techn. Univ.-Sofia, 61, 1, (2011), 9 -- 18.

[13]  I. Trendafilov and D. Vladeva, Endomorphism semirings without zero of a finite chain, Proc. Techn. Univ.-Sofia, 61, 2 (2011) 9 -- 18.

[14]  I. Trendafilov and  D. Vladeva, On some semigroups of the partial transformation semigroup, Appl. Math. in Eng. and Econ. -- 38th Int. Conf. (2012) AIP Conf. Proc. (to appear).

[15]  J. Zumbr\"{a}gel, Classification of finite congruence-simple semirings with zero,
J. Algebra Appl. 7 (2008) 363 -- 377.
}

\vspace{10mm}

Department of Algebra and Geometry, Faculty of Applied Mathematics and Informatics, Technical University of Sofia, 8 Kliment Ohridski Str. Sofia 1000, Bulgaria

\vspace{2mm}

 \emph{e-mail:} ivan$\_$d$\_$trendafilov@abv.bg

\end{document}